# Solving QVIs for Image Restoration with Adaptive Constraint Sets

F. Lenzen[†], J. Lellmann[‡], F. Becker[†], and C. Schnörr[†]

**Abstract.** We consider a class of quasi-variational inequalities (QVIs) for adaptive image restoration, where the adaptivity is described via *solution-dependent* constraint sets. In previous work we studied both theoretical and numerical issues. While we were able to show the existence of solutions for a relatively broad class of problems, we encountered problems concerning uniqueness of the solution as well as convergence of existing algorithms for solving QVIs. In particular, it seemed that with increasing image size the growing condition number of the involved differential operator poses severe problems. In the present paper we prove uniqueness for a larger class of problems and in particular independent of the image size. Moreover, we provide a numerical algorithm with proved convergence. Experimental results support our theoretical findings.

## 1. Introduction.

When solving inverse problems in *image processing* in a variational framework, one faces the issue of selecting a regularizer, which on the one hand should provide suitable reconstruction quality and on the other hand should have sufficient theoretical properties to guarantee existence and uniqueness of a solution.

A common choice is to rely on *convex* regularizers, which in combination with a convex data fidelity term has the advantage that theory and numerics of convex optimization have been intensively studied in literature and are well understood [33]. A prominent example for a convex regularizer with suitable theoretical properties is the total variation (TV) semi-norm [32, 33].

On the other hand, acknowledging the fact that we are reconstructing images, an investigation of the empirical distribution of typical images (Zhu & Mumford [43]) shows, that *non-convex* regularization terms are more appropriate to choose. Using such non-convex regularizers comes with the challenge to prove existence and uniqueness for the resulting variational problem.

In our work ([19, 23, 24] and the present paper), we follow a strategy which combines elements of convex and non-convex formulations. We start with a convex problem using discrete TV regularization and consider its dual formulation, which is given as a

[†]HCI & IPA, Heidelberg University
Speyerer Str. 6, 69115 Heidelberg, Germany
frank.lenzen@iwr.uni-heidelberg.de, becker,schnoerr@math.uni-heidelberg.de
[‡]DAMTP/CIA, University of Cambridge
Centre for Mathematical Sciences, Wilberforce Rd, Cambridge CB3 0WA, United Kingdom
j.lellmann@damtp.cam.ac.uk
The work of J.L. was supported by Award No. KUK-I1-007-43, made by King Abdullah University of Science and Technology (KAUST), EPSRC first grant No. EP/J009539/1, and Royal Society International Exchange Award No. IE110314.

constrained quadratic optimization problem

$$\min_{p \in \mathcal{D}} F(p) \tag{1.1}$$

with a convex constraint set $\mathcal{D}$. The next step is to make the constraint set $\mathcal{D}$ depending on a fixed element, say $p_0$, i.e. we consider $\mathcal{D} = \mathcal{D}(p_0)$. This generalization makes the regularization approach an adaptive one. Since the optimization problem is still convex, existence of a solution $\overline{p}$ is guaranteed. Finally, we consider the problem of finding a fixed-point $p^*$ of the mapping

$$p_0 \to \overline{p} := \arg\min_{p \in \mathcal{D}(p_0)} F(p). \tag{1.2}$$

We refer to the convex problem in (1.2), where $p_0$ is fixed, as the *inner* problem.

The *overall problem* of finding a fixed-point $p^*$ is equivalent to solving a quasi-variational inequality (QVI) [10], thus we can make use of available theory on existence. Uniqueness of such a fixed-point, however, in general is an open issue, when the underlying dual functional ($p_0$ fixed) is not *strictly* convex. To tackle this uniqueness issue under suitable conditions is one of the main contributions of this paper. In particular, it turns out that under the assumption that the mapping $p_0 \to \overline{p}$ is a contraction w.r.t. a suitable semi-norm, the corresponding primal solution is unique. Having found the fixed-point $p^*$, this fixed-point defines the adaptivity of the regularizer via the constraint set $\mathcal{D}(p^*)$, while also being the solution of the convex inner problem. As a consequence, our approach implements a solution-driven adaptivity instead of a data-driven one. Furthermore, we can investigate the behavior of the regularizer at the given fixed-point and find that it mimics the behavior of a non-convex regularizer.

**Related work:** As exemplary applications for our regularization approach we consider total variation based image denoising and non-blind image deblurring. We start with related work concerning image denoising. For the task of denoising total variation regularization was introduced by Rudin, Osher and Fatemi in [32]. Various modifications have been proposed to make this functional adaptive to the input data [2, 4, 14, 18, 34, 36, 42]. The approaches in [2, 4, 14, 36] can be described by means of locally dependent constraint sets (data-driven), i.e. a fixed $p_0$ in our formulation.

Another important class of approaches are the non-local methods [7, 17, 20, 31] including non-local variants of TV. These non-local variants can also be regarded as adaptive, since their local weights are depending on the input. On the other hand, adaptive methods which steer adaptivity by locally averaging the input data over a neighborhood, e.g. using the structure tensor [16], can be interpreted as non-local.

Another commonly used modification of the ROF functional is to replace the $L_2$ norm of the data fidelity term by an $L_1$ norm [1]. We remark that, by using a standard splitting of variables, the approach presented here can also be formulated with such an $L_1$ data term.

Recent developments in the field of TV regularization focus also on extending TV to second- or higher-order [6, 35]. In [22] we have proposed an anisotropic approach of

first- and second-order TV, which due to its formulation by varying constraint sets also fits into the concept of solution-driven adaptivity presented here.

The regularization approaches considered above can also be used for the task of deblurring, see e.g. [8, 11, 12]. In [29] an TV based deblurring approach with adaptive choice of the regularization parameter has been proposed. Similar to the task of denoising, non-local operators have being also considered for TV deblurring, see e.g. [20, 38].

Besides for image restoration tasks, TV-based regularization approaches are widely used for other inverse problems in computer vision, e.g. for optical flow [41, 40] and multi-labeling [21, 37].

Concerning image restoration with non-convex regularization, in addition to Zhu & Mumford [43] we want to mention here the work by Charbonnier and co-workers [13, 5] and by Levin [25].

As already mentioned, the fixed-point problem (1.2), which is the core problem in our considerations, is equivalent to a quasi-variational inequality. We make use of the work on QVIs presented in [10, 27, 28]. While theory on existence can be directly utilized, uniqueness results do not apply due to the non-strict convexity of the inner problem. We discuss these issues in detail in the main part of this paper.

**Contribution:** In our previous work [23, 24] we sketched the proposed framework and provided existence theory. In [24] we showed uniqueness for a very narrow class of problems, which also scaled unfavorably with the image size.

In the present paper, we show uniqueness of a fixed-point for a broad class of QVIs for image denoising, namely those, for which the underlying solution operator is a contraction. In particular, this condition is not depending on the image size. Thus, our theoretical results significantly generalize our previous work in [24].

Moreover, we give a detailed discussion why classical results (Noor et al. [28] and Nesterov & Scrimali [27]) on the uniqueness of solutions of QVIs can not directly be applied to our framework. However, there is a strong relationship between our theoretical considerations and the work of [27, 28].

Finally, we propose an algorithm for solving the considered QVIs and prove convergence. We support our theoretical results by numerical experiments.

**Paper organization:** Our paper is organized as follows. We start with a review of three case examples of TV regularization for image denoising and non-blind image deblurring in Sect. 2. In Sect. 3 we recall our model of solution-driven adaptivity described by means of quasi variational inequalities. We consider theoretical results in Sect. 4, where we firstly recall theory on existence (Sect. 4.1), then discuss the impact of existing work on uniqueness (Sect. 4.2.1) and finally prove uniqueness for the considered QVIs under suitable conditions (Sect. 4.2.2). In Sect. 5 we provide an algorithm and prove its convergence. We present numerical experiments supporting our theoretical results in Sect. 6.

**2. TV-Regularization and Data-Driven Adaptivity.** In the following, we recall several variational approaches for image denoising and non-blind image deblurring that are based on total variation (TV) regularization. These approaches will be the starting point for our generalizations in Section 3.

We use the following general notations. Firstly, let $\Omega \subset \mathbb{R}^d$ be a $d$-dimensional open, bounded domain with Lipschitz boundary. Secondly, in $\mathbb{R}^n$, $n$ arbitrary, we denote the closed ball with radius $\alpha$ centered at 0 by $B_\alpha(0)$.

**2.1. Image Denoising.** We consider the standard noise model, where some noise-free image $u$ is distorted by additive i.i.d. Gaussian noise with zero mean. For the noisy image we use the notation $f$ and assume $f \in L^2(\Omega)$. We refer to $u$ as the *original* or *ground truth* image.

**2.1.1. The Classical ROF Model.** We start with the total variation denoising approach by Rudin, Osher & Fatemi (ROF) [32],

$$\min_{u \in BV(\Omega)} E(u), \qquad E(u) := \frac{1}{2}\|u - f\|_{L^2}^2 + \alpha \operatorname{TV}(u), \qquad \alpha > 0, \tag{2.1}$$

where $BV(\Omega), \Omega \subset \mathbb{R}^d$ is the space of functions of bounded total variation and

$$\operatorname{TV}(u) := \sup \left\{ \int_\Omega u(x) \operatorname{div} p(x)\, dx \mid p \in C_c^\infty(\Omega; \mathbb{R}^d) \colon \|p(x)\|_2 \leq 1 \right\}. \tag{2.2}$$

is the total variation semi-norm. We rewrite $\alpha \operatorname{TV}(u)$ in terms of *constraint sets*:

$$\begin{aligned} \alpha \operatorname{TV}(u) &= \sup\{\langle u, v\rangle_{L^2} \mid v \in \mathcal{C}\}, \\ \mathcal{C} &:= \operatorname{div} \mathcal{D}, \\ \mathcal{D} &:= \{p \in C_c^\infty(\Omega; \mathbb{R}^d) \colon \|p(x)\|_2 \leq \alpha\}, \end{aligned} \tag{2.3}$$

where div is applied element-wise on $\mathcal{D}$. The dual problem of (2.1) (cf. [8]) can be formulated as

$$\min_{p \in \overline{\mathcal{D}}} \frac{1}{2}\|f - \operatorname{div} p\|_{L^2}^2, \tag{2.4}$$

where $\overline{\mathcal{D}}$ is the closure of $\mathcal{D}$.

Let us now consider a discretization of (2.1). To this end, we consider an equidistant grid on $\Omega$ with $n$ grid points. The grid values of the dual variable $p(x) \in \mathbb{R}^d$ are interpreted as a vector $p \in \mathbb{R}^{nd}$. The dual problem (2.4) in the discrete formulation then becomes

$$\min_{p \in \mathcal{D}} F(p), \qquad F(p) := \frac{1}{2}\|\operatorname{L} p - f\|_2^2, \tag{2.5}$$

where $\operatorname{L} : \mathbb{R}^{nd} \to \mathbb{R}^n$ is a discretization of the divergence operator div. The constraint set in (2.4) becomes

$$\mathcal{D} = \mathcal{D}_1^{loc} \times \mathcal{D}_2^{loc} \times \cdots \times \mathcal{D}_n^{loc}, \tag{2.6}$$

where each *local constraint set* $\mathcal{D}_i^{loc}$, $i = 1, \ldots, n$ is a $d$-dimensional closed ball $B_\alpha(0)$ of radius $\alpha$. The dual problem (2.5) will be the starting point for our generalization in Sect. 3.

**2.1.2. Higher Order Total Variation.** Analogously to TV regularization of first-order, higher-order models can be considered. We exemplarily focus on the second order total variation in the case $d = 2$:

$$\mathrm{TV}^2(u) := \sup\left\{ \int_\Omega u(x)\,\mathrm{div}^2\, p(x)\; dx \mid p \in C_c^\infty(\Omega;\mathbb{R}^4) \colon \|p(x)\|_2 \le 1 \right\}, \tag{2.7}$$

where $\mathrm{div}^2\, p := \partial_{xx} p_1 + \partial_{xy} p_2 + \partial_{yx} p_3 + \partial_{yy} p_4$. Typically first- and second-order TV are used jointly for regularization, e.g., for the task of denoising, one could solve

$$\min_{u \in BV^2(\Omega)} E(u), \qquad E(u) := \|u - f\|_{L^2}^2 + \alpha\,\mathrm{TV}(u) + \beta\,\mathrm{TV}^2(u), \qquad \alpha > 0, \beta > 0, \tag{2.8}$$

where $BV^2(\Omega)$ is the space of functions with bounded total variation of first- and second-order (see [33, Section 9.8] for details).

Proceeding analogously to the case of first-order TV, we can derive a dual formulation of (2.8), which after discretization reads similar to (2.5):

$$\min_{p \in \mathcal{D}} F(p), \qquad F(p) := \frac{1}{2}\|\,\mathrm{A}\,p - f\|_2^2, \tag{2.9}$$

where $\mathrm{A} : \mathbb{R}^{6n} \to \mathbb{R}^n$ is a operator discretizing $\mathrm{div}\, p_1 + \mathrm{div}^2\, p_2$ with $p = (p_1, p_2) \in \mathbb{R}^{2n} \times \mathbb{R}^{4n}$. The constraint set $\mathcal{D}$ in (2.9) is given by a product set of local constraint sets $\mathcal{D}_i^{loc}$, where each set $\mathcal{D}_i^{loc}$ is again a product of a two-dimensional ball $B_\alpha(0)$ of radius $\alpha$ and a four-dimensional ball $B_\beta(0)$ of radius $\beta$.

We refer to [6] for the alternative model of Total Generalized Variation (TGV), which is based on a different operator A and a different constraint set $\mathcal{D}$.

**2.2. Image Deblurring.** In this section we consider the task of image deblurring/deconvolution. For the sake of simplicity, we focus on non-blind deconvolution, where the convolution kernel is known a-priori. The problem formulation is as follows. Let $f$ be some observed data, which are obtained from a noise-free image $u$ by convolution with a kernel $M(x) : \Omega \to \mathbb{R}$, followed by an addition of Gaussian noise, i.e.

$$f = M * u + \delta, \tag{2.10}$$

where $\delta$ is a realization of a Gaussian random variable with zero mean.

In order to recover $u$ from $f$, assuming that $u \mapsto M * u$ is an operator mapping from $L^2(\Omega) \to L^2(\Omega)$, we aim at minimizing

$$\operatorname*{arg\,min}_{u \in L^2(\Omega) \cap BV(\Omega)} \; \tfrac{1}{2}\|M * u - f\|_{L^2}^2 + \alpha\,\mathrm{TV}(u). \tag{2.11}$$

Moving to a discrete formulation of the problem, we now assume that $u, f \in \mathbb{R}^n$ are the function values at the $n$ nodes of an equidistant two-dimensional grid. Moreover, we replace the continuous convolution $M * u$ by a matrix-vector-product $Mu$, where $M$ now denotes a $n \times n$ matrix. In what follows we assume that $M$ is invertible. As in the previous examples, we denote by L the discretization of the divergence operator div. The optimization problem we consider is given as

$$\arg\min_{u\in\mathbb{R}^n} E(u) = \tfrac{1}{2}\|Mu - f\|_2^2 + \sup_{p\in\mathcal{D}} (\mathrm{L}\,p)^\top u, \tag{2.12}$$

where

$$\mathcal{D} = \{p \in \mathbb{R}^{nd}, p_i \in B_\alpha(0), i = 1,\dots,n\} \tag{2.13}$$

for $p = (p_1, p_2, \dots, p_n)^\top$ with $p_i \in \mathbb{R}^2$. We derive the corresponding dual problem as follows. The optimality condition for $u$ reads

$$M^\top(Mu - f) + \mathrm{L}\,p = 0. \tag{2.14}$$

We deduce from (2.14) that

$$\begin{aligned} u &= (M^{-1}M^{-\top})(M^\top f - \mathrm{L}\,p) = M^{-1}(f - M^{-\top}\,\mathrm{L}\,p),\\ Mu &= f - M^{-\top}\,\mathrm{L}\,p, \end{aligned} \tag{2.15}$$

where $M^{-\top} := (M^\top)^{-1}$. Inserting (2.15) in (2.12) and using the abbreviation $\mathrm{A} := (M^{-\top}\,\mathrm{L})$, we obtain

$$\begin{aligned} E^*(p) ={}& \tfrac{1}{2}\|Mu - f\|_2^2 + (\mathrm{L}\,p)^\top u && (2.16)\\ ={}& \tfrac{1}{2}\|M^{-\top}\,\mathrm{L}\,p\|_2^2 + (M^{-1}f)^\top\,\mathrm{L}\,p - (M^{-1}M^{-\top}\,\mathrm{L}\,p)^\top\,\mathrm{L}\,p && (2.17)\\ ={}& \tfrac{1}{2}\|M^{-\top}\,\mathrm{L}\,p\|_2^2 + f^\top M^{-\top}\,\mathrm{L}\,p - (M^{-1}M^{-\top}\,\mathrm{L}\,p)^\top\,\mathrm{L}\,p && (2.18)\\ ={}& -\tfrac{1}{2}\|\mathrm{A}\,p\|_2^2 + f\,\mathrm{A}\,p - (\mathrm{A}\,p)^\top\,\mathrm{A}\,p && (2.19)\\ ={}& -\tfrac{1}{2}\|\mathrm{A}\,p - f\|_2^2 + \tfrac{1}{2}\|f\|_2^2. && (2.20) \end{aligned}$$

When maximizing $E^*(p)$ over $\mathcal{D} = \{p \in \mathbb{R}^{nd}, p_i \in B_\alpha(0)\}$, the constant term $\frac{1}{2}\|f\|_2^2$ can be omitted without changing the optimum. Moreover, switching from the maximization of $E^*$ to the minimization of $F(p) := -E^*(p)$, we can formulate the dual problem of (2.12) as

$$\arg\min_{p\in\mathcal{D}} F(p), \qquad F(p) = \tfrac{1}{2}\|\mathrm{A}\,p - f\|_2^2. \tag{2.21}$$

From a solution $\overline{p}$ of the dual problem we can retrieve the solution $\overline{u}$ of the primal problem by by $\overline{u} = M^{-1}(f - \mathrm{A}\,\overline{p})$.

We observe that the dual problem attains the same form as in the examples before (cf. Eqns. (2.5) and (2.9)).

**2.3. Adaptive Regularization.** In the literature various adaptive TV approaches have been proposed. They can generally be divided into two classes, namely, approaches with *locally varying* regularization strength and *anisotropic* TV approaches. Both concepts are covered by the formulation via constraint sets as follows. Starting with the general form

$$\arg\min_{p\in\mathcal{D}} F(p), \qquad F(p) = \tfrac{1}{2}\|\mathrm{A}\,p - f\|_2^2, \tag{2.22}$$

where $\mathcal{D}$ is the product set of the local constraint sets $\mathcal{D}_1^{loc}, \dots, \mathcal{D}_n^{loc}$ as in (2.6), we now allow the sets $\mathcal{D}_i^{loc}$ to vary locally:

- By individually changing the *size* of $\mathcal{D}_i^{loc}$, e.g. depending on the noise or image content, the regularization strength changes locally.
- By choosing *anisotropic shapes* for $\mathcal{D}_i^{loc}$, e.g. rectangles [4], parallelograms [36], and ellipses [2, 19, 22], a directionally dependent regularization is introduced.

In both cases, the introduced adaptivity has to be steered by additional information, e.g. about noise level, edge position and edge orientation. The standard way is to either estimate the required properties as additional unknowns in the optimization process, or to examining a pre-smoothed version of the data $f$. The later case formally can be regarded as introducing a dependency of $\mathcal{D}$ on $f$, i.e. $\mathcal{D} = \mathcal{D}(f)$.

We pick up three different examples of adaptive/anisotropic TV regularization which follow the latter concept. The first is obtained from the standard ROF model by locally varying parameter $\alpha$.

Example 2.1 (Data-driven adaptivity). *Let us consider the following generalization of the optimization problems* (2.1) *and* (2.11)*, where the parameter $\alpha$ is allowed to change locally:*

$$E(u) = \int_\Omega \frac{1}{2}(M * u - f)^2 \, dx + \sup\{\int_\Omega u \operatorname{div} p \, dx \mid p \in C_c^1(\Omega, \mathbb{R}^d), |p(x)| \leq \alpha(x)\}, \quad (2.23)$$

*where we assume $\alpha(x) \geq c > 0$. Our aim is to reduce the local regularization parameter $\alpha(x)$ at edges. A simple way to find such edges would be to consider the gradient magnitude of the input data $f$ and to set $\alpha(x) := \max(\alpha_0(1 - \kappa|\nabla f(x)|), \varepsilon)$ with a constant $\alpha_0$ determining the maximal regularization strength and some small $\varepsilon > 0$ to ensure boundedness of $\alpha$ from below by a positive constant, which ensures existence [22].*

*However, since the approach should be robust against the noise contained in $f$, a pre-smoothing of $f$ before evaluating the gradient is inevitable.To this end, let $f_\sigma := K_\sigma * f$ be the convolution of $f$with a Gaussian kernel $K_\sigma$ with standard deviation $\sigma > 0$. An adaptive choice of $\alpha(x)$ is*

$$\alpha(x) := \max\{\alpha_0(1 - \kappa|(\nabla f_\sigma(x))|), \varepsilon\}. \quad (2.24)$$

*There exist alternative choices for varying the regularization strength, such as the g-function from the Perona-Malik model [30], or models utilizing the structure tensor [16].*

*Considering again the dual formulation of* (2.23) *in a discrete setting, we retain the form* (2.5)*, with the only difference that the local constraint sets become dependent on the spatial location and the input data $f$,*

$$\mathcal{D}_i^{loc} := B_{\alpha_i}(0), \qquad \alpha_i := \max\{\alpha_0(1 - \kappa|(\mathrm{L}^\top f_\sigma)_i|), \varepsilon\}, \quad (2.25)$$

*where $i = 1, \dots, n$ are the indices of the grid nodes. Note that in the discrete setting, where $\mathrm{L}$ is a discretization of* div*, the discrete pendant of $\nabla$ is $-\mathrm{L}^\top$. Recall that the dual problem reads*

$$\arg\min_{p \in \mathcal{D}} \tfrac{1}{2}\| \mathrm{A}\, u - f\|_2^2 \quad (2.26)$$

*with*

$$\mathcal{D} := \mathcal{D}_1^{loc} \times \cdots \times \mathcal{D}_n^{loc} \quad (2.27)$$

*and* $\mathrm{A} = \mathrm{L}$ *in the case of denoising and* $\mathrm{A} = (M^{-\top})\,\mathrm{L}$ *in the case of deblurring.*

*The key observation in this example is that we formally introduced a dependency of* $\mathcal{D}$ *on* $f$ *via* $\alpha$. *We denote this dependency by* $\mathcal{D}(f)$. *We refer to this concept as* data-driven adaptivity. ♦

Our second example generalizes the ROF model by considering a directionally dependent regularization, which results in an anisotropic shape of the local constraint sets.

**Example 2.2 (Anisotropic first-order TV).** *We consider an anisotropic TV regularization with a strong penalization of the image gradient in homogeneous regions (isotropic) and, at edges, a weak penalization in normal direction and a strong penalization in tangential direction to the edge (anisotropic).*

*To this end we require information about the location and orientation of edges in terms of an edge indicator function* $\chi_e : \Omega \to [0,1]$ *and a vector field* $v_e : \Omega \to \mathbb{R}^2$ *of edge normals, which both can be obtained from the standard structure tensor [16] of* $f$ *by setting*

$$\chi_e(x; f) = \min\{\kappa(\lambda_1(x) - \lambda_2(x)), 1\}, \tag{2.28}$$

$$v_e(x; f) = w_1(x), \tag{2.29}$$

*where* $\lambda_1 \geq \lambda_2 \geq 0$ *are the ordered eigenvalues of the structure tensor,* $w_1$ *is the eigenvector to eigenvalue* $\lambda_1$ *and* $\kappa > 0$ *is a parameter controlling the edge sensitivity. We refer to [24] for exact definitions and further details.*

*With this edge information, we choose* $\mathcal{D}_i^{loc} = \mathcal{D}^{loc}(x_i)$ *at grid node* $x_i$ *to be an ellipse with one half axis parallel to* $v_e(x_i, f)$ *of length* $\chi_e(x_i)\alpha + (1 - \chi_e(x_i))\beta$, *with constants* $0 \leq \alpha \leq \beta$, *and the perpendicular half axis of length* $\beta$.

*The cross product of the local constraint sets* $\mathcal{D}_i^{loc}$ *as in (2.27) defines our (global) constraint set* $\mathcal{D}(f)$. ♦

Finally, let us consider an example of adaptive higher-order TV regularization.

**Example 2.3 (Adaptive first- and second-order TV).** *We revisit the first- and second-order TV regularization models from Sect. 2.1.2 with the discretized dual problem*

$$\min_{p \in \mathcal{D}} F(p), \qquad F(p) := \frac{1}{2}\|\,\mathrm{A}\,p - f\|_2^2, \tag{2.30}$$

*where the operator* $\mathrm{A} : \mathbb{R}^{6n} \to \mathbb{R}^n$ *discretizes* $\operatorname{div} p_1 + \operatorname{div}^2 p_2$ *with* $p = (p_1, p_2) \in \mathbb{R}^{2n} \times \mathbb{R}^{4n}$.

*We are aiming at a regularization with locally varying regularization strengths* $\alpha_i$ *for first- and* $\beta_i$ *for second-order. Analogously to Example 2.1, we choose*

$$\alpha_i := \max\{\alpha_0(1 - \kappa|(\mathrm{L}^\top f_\sigma)_i|), \varepsilon\}, \tag{2.31}$$

$$\beta_i := \max\{\beta_0(1 - \kappa|(\mathrm{L}^\top f_\sigma)_i|), \varepsilon\}, \tag{2.32}$$

*with constants* $\alpha_0, \beta_0 > 0$, *i.e. in homogeneous regions (vanishing gradient* $\nabla f = 0$*) we penalize the first- and second-order TV with factor* $\alpha_0$ *and* $\beta_0$, *respectively, while we* *reduce the regularization strength at edges ($|\nabla f| \gg 0$). As local constraint sets we then choose*

$$\mathcal{D}_i^{loc} := B_{\alpha_i}(0) \times B_{\beta_i}(0) \subset \mathbb{R}^2 \times \mathbb{R}^4. \tag{2.33}$$

♦

The above examples show, that many popular variational approaches conform to the generic model (2.22). A limitation of the above adaptive approaches is, that the adaptivity is determined by the noisy input data $f$ (data-driven adaptivity), rather than by the noise-free solution $u$. In the next section, we show how we can switch from a data-driven to a solution-driven adaptivity.

## 3. Solution-driven Adaptivity.

In [24] we have proposed a new kind of adaptivity, where the constraint set $\mathcal{D}$ depends on the unknown solution of the problem. We recall this approach below.

Our approach generalizes the examples of Section 2 with respect to the operator $\mathrm{A}$ and the form of the constraint set $\mathcal{D}$.

We describe this generalization in a discrete setting, where we consider again an equi-distant grid on $\Omega$ with $n$ grid points. We start with a dual problem of the form

$$\min_{p \in \mathcal{D}} F(p), \qquad F(p) := \frac{1}{2} \| \mathrm{A}\, p - f\|_2^2, \tag{3.1}$$

where $\mathrm{A} : \mathbb{R}^{mn} \to \mathbb{R}^n$ now is a general discrete operator. We assume that the constraint set $\mathcal{D}$ takes the form

$$\mathcal{D} = \mathcal{D}_1^{loc} \times \mathcal{D}_2^{loc} \times \cdots \times \mathcal{D}_n^{loc}, \tag{3.2}$$

where each $\mathcal{D}_i^{loc}$ is a local $m$-dimensional closed convex constraint set at the $i$-th grid point. We stress that the shape of $\mathcal{D}^{loc}$ can be arbitrary. The solution of the primal problem can be retrieved by $\overline{u} := M^{-1}(f - \mathrm{A}\,\overline{p})$ from the solution $\overline{p}$ of the dual problem (3.1).

We remark that the dual problem (3.1) can be equivalently formulated based on a variational inequality (VI)

$$\text{find} \quad \overline{p} \in \mathcal{D} \text{ such that } \langle \nabla F(\overline{p}), p - \overline{p} \rangle \geq 0 \quad \forall p \in \mathcal{D}. \tag{3.3}$$

In our case, the gradient of $F(p)$ is an affine function of $p$:

$$\nabla F(p) = \mathrm{A}^\top (\mathrm{A}\, p - f). \tag{3.4}$$

We will make use of this specific form in the following section.

We now generalize the problem (2.5) by introducing a dependency of $\mathcal{D}$ on the dual variable, i.e. $\mathcal{D} = \mathcal{D}(p_0)$ for some $p_0 \in \mathbb{R}^{mn}$ and search for a fixed-point $p^*$ of the mapping

$$\boxed{p_0 \mapsto \overline{p} := \arg\min_{p \in \mathcal{D}(p_0)} F(p).} \tag{3.5}$$

Please note that we have to distinguish between a fixed-point of (3.5), denoted by $p^*$, and a minimizer $\overline{p}$ of the convex dual problem $\arg\min_{p\in\mathcal{D}(p_0)} F(p)$ for fixed $p_0$. Both coincide only if $p_0 = p^*$.

Having found a fixed point $p^*$, the corresponding constraint set is $\mathcal{D}(p^*)$, i.e. the adaptivity becomes *solution-driven*.

Moreover, we can interpret $p^*$ as the solution of a convex problem with fixed constraint set $\mathcal{D} = \mathcal{D}(p^*)$ and can consider the solution $u^*$ of the corresponding primal problem, which can be retrieved by $u^* = M^{-1}(f - \mathrm{A}\,p^*)$ (with $M := \mathrm{Id}$ in the case of denoising).

Introducing the fixed-point problem (3.5) has several advantages:

1. The *inner problem*, i.e. the problem of finding $\arg\min_{p\in\mathcal{D}(p_0)} F(p)$ for a fixed $p_0$ is a convex problem. Theoretical and numerical issues of this problem have been intensively studied.
2. Also for the outer fixed-point problem, theory on existence is at hand.
3. Concerning the inner problem, our ansatz allows us to switch between primal, dual and the saddle-point formulation

$$\min_{u\in\mathbb{R}^n}\left\{\frac{1}{2}\|Mu - f\|_2^2 + \sup_{p\in\mathcal{D}(p_0)} \langle u, \mathrm{L}\,p\rangle_{L^2}\right\}, \tag{3.6}$$

for fixed $p_0$. In particular, after having found the fixed-point $p^*$, we can retrieve the primal solution $u^*$ as the solution of (3.6) with fixed $p_0 = p^*$.

Remark 3.1. *The concept of a solution-driven adaptivity also covers the case that the adaptivity is determined based on the primal variable $u$, since we can express $\mathcal{D}(u)$ by $\mathcal{D}(p)$ using the relationship $u = M^{-1}(f - \mathrm{A}\,p)$. However, the fixed-point problem (3.5) in general is not equivalent to the non-convex problem*

$$\arg\min_{u\in\mathbb{R}^n} \|Mu - f\|_{L^2}^2 + \sup_{p\in\mathcal{D}(u)} (\mathrm{L}\,p)^\top u. \tag{3.7}$$

Let us illustrate the considerations made so far by an example:

Example 3.2. *We compare the two conceptually different ways of implementing adaptivity – data-driven adaptivity, where $\mathcal{D}$ depends solely on the input data $f$, and solution-driven adaptivity, where the constraint set $\mathcal{D}$ depends on the unknown $u$ (or, equivalently, $p$).*

*Firstly, we recall the data-driven adaptive TV regularization from Example 2.1, where the regularization parameter $\alpha$ was chosen locally at grid node $i$ to be*

$$\alpha_i := \max\{\alpha_0(1 - \kappa|(\mathrm{L}^\top f_\sigma)_i|), \varepsilon\}. \tag{3.8}$$

*Our proposed generalized approach permits to make the constraint set depending on $u$. To this end, let*

$$\alpha_i(u) := \max\{\alpha_0(1 - \kappa|(\mathrm{L}^\top u)_i|), \varepsilon\} \tag{3.9}$$

*(assuming that $u$ is noise-free, we omit the Gaussian pre-smoothing), and*

$$\mathcal{D}(u) = \mathcal{D}_1^{loc} \times \mathcal{D}_2^{loc} \times \cdots \times \mathcal{D}_n^{loc}, \text{ where } \mathcal{D}_i^{loc}(u) := B_{\alpha_i(u)}(0). \tag{3.10}$$

*Considering alternatively the dual problem (2.5), we can by means of the relationship $u = M^{-\top}(f - \mathrm{A}\,p)$ instead assume that $\mathcal{D}$ depends on $p$, or, more precisely, on $\mathrm{A}\,p$:*

$$\mathcal{D} = \mathcal{D}(p) = \mathcal{D}(\mathrm{A}\,p). \tag{3.11}$$

*Although there is also a formal dependency on $f$, we omit this in our notation to emphasize the different models $\mathcal{D}(f)$ (adaptive to data $f$, but fixed) and $\mathcal{D}(p)$ (adaptive to the unknown $p$).*

*We will compare both models experimentally in Section 6.* ♦

As already mentioned before, alternative choices for varying the regularization strength $\alpha$, such as using the function $g(|\nabla u|)$ from the Perona-Malik diffusion model [30], exist. In view of the theory provided in the next section, such an $\alpha(u)$ should at least be Lipschitz-continuous w.r.t. $u$.

Moreover, we stress that besides the examples discussed in Sect. 2.3 various other models of adaptive/anisotropic regularization exists, which are covered by the above general model (3.5), see e.g. [4, 22, 23, 24, 36].

Finally, we remark that the generalized problem of finding a fixed-point of (3.5) is equivalent to solving a quasi-variational inequality problem (QVIP) (cf. [10])

$$\boxed{\text{find } p^* \in \mathcal{D}(p^*) \text{ such that } \langle \nabla F(p^*), p - p^* \rangle \geq 0 \quad \forall p \in \mathcal{D}(p^*)} \tag{3.12}$$

with $F(p) = \frac{1}{2}\|\mathrm{A}\,p - f\|^2$. When reformulating the proposed fixed-point problem as a QVIP, we can make use of the theory existing in literature [10, 27].

We will provide existence and uniqueness results for the QVIP (3.12) in detail in the subsequent section.

## 4. Theory.

The key issue of this section is to prove uniqueness for the problem (3.12) under sufficient conditions. A prerequisite for uniqueness is the existence of a solution. We therefore briefly recall existence results from literature in the next section, before turning to uniqueness results in Section 4.2.

### 4.1. Existence.

We recall existence results from [24] for problem (3.12) together with the necessary assumptions. These assumptions will also be required for uniqueness results provided in second part of this section.

Assumption 4.1 (for existence). *Assume that*

$$\mathcal{D}\colon \overline{p} \quad \rightrightarrows \quad \mathcal{D}(\overline{p}) := \{p \in \mathbb{R}^{mn} \colon p_i \in \mathcal{D}_{loc}^i(\overline{p}) \subset \mathbb{R}^{mn},\ i = 1, \ldots, n\}, \tag{4.1}$$

*where each $\mathcal{D}_{loc}^i : \mathbb{R}^{mn} \rightrightarrows \mathbb{R}^{mn}, i = 1, \ldots, n$ has the following properties:*

*(i) For fixed $p$ the set $\mathcal{D}_{loc}^i(p)$ is a closed convex subset of $\mathbb{R}^{mn}$.*

*(ii)* *There exists $C > 0$, such that for all $i, p$: $\mathcal{D}^i_{loc}(p) \subset B_C(0)$ (closed ball with radius $C$).*
*(iii)* *There exists $c > 0$, such that for every $p$ and every $i$ we have $B_c(0) \subset \mathcal{D}^i_{loc}(p)$. In particular, $\mathcal{D}^i_{loc}(p)$ is non-empty.*
*(iv)* *The projection $\Pi_{\mathcal{D}^i_{loc}(p)}(q)$ of $q$ onto $\mathcal{D}^i_{loc}(p)$ for a fixed $q$ is continuous w.r.t. $p$.*

Proposition 4.2. *Let $F(p) := \frac{1}{2}\|f - \mathrm{A}\, p\|_2^2$, where $\mathrm{A} : \mathbb{R}^{mn} \to \mathbb{R}^n$ is a linear operator. Moreover, let $\mathcal{D}(p)$ be defined as in (4.1), such that $\mathcal{D}^i_{loc}(p), i = 1, \ldots, n$ satisfy Assumption 4.1. Then the problem (3.12) has a solution.*

*Proof.* See [24, Prop. 1]. ■

The proof in [24] utilizes a general existence result for QVIs presented in [10], whose core ingredient is Brouwer's fixed-point theorem and which makes use of the continuity of the mapping $p \to \Pi_{\mathcal{D}(p)}(q)$ (guaranteed by Assumption 4.1(iii)). We will see that for uniqueness results, a higher regularity of $p \to \Pi_{\mathcal{D}(p)}(q)$, namely a Lipschitz-continuity is required.

Remark 4.3 (A-priori bounds). *From Assumption 4.1 (ii) we derive an a-priori bound for $\mathcal{D}(p)$ independent from $p$:*

$$\mathcal{D}(p) \subset (B_C(0))^n \subset B_{\sqrt{n}C}(0) \qquad \forall p \in \mathbb{R}^{mn}. \tag{4.2}$$

*We define $R := \sqrt{n}C$. In particular, (4.2) provides a bound for a solution $p^*$ of (3.12):*

$$p^* \in \mathcal{D}(p^*) \subset B_R(0), \;\; i.e.\;\; \|p^*\|_2 \leq R = \sqrt{n}C. \tag{4.3}$$

**4.2. Uniqueness.** Let us now consider uniqueness results for the QVI (3.12). This part comprises the main contribution of this paper.

We start with a discussion on related work in Sect. 4.2.1, in particular the paper by Nesterov & Scrimali [27], which provides existence results for strongly monotone gradients $\nabla F$ under certain conditions. We will see that this theory is only partially applicable in our context, since $\nabla F$ in our case is strongly monotone only on a subspace of $\mathbb{R}^{mn}$. Consequently, we will be able to show uniqueness of $p$ only with respect to its component in that subspace. The final uniqueness result is provided in Sect. 4.2.2.

**4.2.1. Existing Theory.** We recall QVI (3.12), which is of the form

$$\langle g(p^*), p - p^* \rangle \geq 0 \qquad \forall p \in \mathcal{D}(p^*), \tag{4.4}$$

with $g(p) = \nabla F(p) = \mathrm{A}^\top(\mathrm{A}\, p - f)$.

Theory on existence of a unique solution $p^* \in \mathcal{D}(p^*)$ to (4.4) under certain conditions has been shown by Noor & Oettli [28, Thm. 9] and Nesterov & Scrimali [27, Cor. 2]. We focus on the latter as their results are more general.

We briefly recall the required conditions below. It is of particular importance that these conditions have to hold for an arbitrary norm $\|x\|_B := \sqrt{x^\top B x}$ for a positive definite matrix $B$. Note that the scalar product in (3.12) is the standard product independent from $B$.

Firstly, operator $g : \mathbb{R}^{mn} \to \mathbb{R}^{mn}$ is assumed to be Lipschitz-continuous with parameter $\mu_B > 0$, i.e.

$$\|g(x) - g(y)\|_{B^*} \leq \mu_B \|x - y\|_B \quad \forall x, y \in \mathbb{R}^{mn}, \tag{4.5}$$

where $\|.\|_{B^*}$ is the norm in the dual space of $\mathbb{R}^{mn}$ equipped with $\|.\|_B$. Note that constant $\mu_B$ depends on the chosen norm $B$. We indicate this dependency by the subscript $B$.

Secondly, $g$ is assumed to be strongly monotone with parameter $\nu_B$, again depending on $B$, i.e.,

$$\langle g(x) - g(y), x - y\rangle \geq \nu_B \|x - y\|_B^2 \quad \forall x, y \in \mathbb{R}^{mn}. \tag{4.6}$$

Both constants $\mu_B$ and $\nu_B$ define the condition number $\gamma_B := \frac{\mu_B}{\nu_B}$, which in our case is the condition number of $\mathrm{A}^\top \mathrm{A}$.

Finally, it is assumed that the projection $\Pi_{\mathcal{D}(p)}(q)$ is Lipschitz-continuous w.r.t. $p$, i.e. for arbitrary $q \in \mathbb{R}^{mn}$,

$$\|\Pi_{\mathcal{D}(p)}(q) - \Pi_{\mathcal{D}(\tilde{p})}(q)\|_B \leq \eta_B \|p - \tilde{p}\|_B. \tag{4.7}$$

We refer to $\eta_B$ as the *variation rate* of $\mathcal{D}(p)$.

In the following, we use the notation $\mu_2$, $\nu_2$ and $\eta_2$ whenever we are considering the standard Euclidean norm $\|\cdot\|_B = \|\cdot\|_2$ ($B = \mathrm{Id}$).

Under the above assumptions, Cor. 2 in [27] provides uniqueness in the case that

$$\eta_B \gamma_B < 1. \tag{4.8}$$

One immediately observes that two open issues preclude the direct application of the theory in [27, 28] to our problem:

- Operator $\nabla F$ in (3.12) has a non-trivial null space $\mathcal{N}(\mathrm{A})$ and thus is not strongly monotone.
- On the complement $\mathcal{N}^\perp(\mathrm{A})$ of the null space, the condition number $\gamma_2$ w.r.t. the standard Euclidean norm tends to infinity with increasing problem size. As a consequence, assuming that $\eta_2$ is fixed, (4.8) can not be satisfied for arbitrary large image. Alternatively, in order to guarantee (4.8), $\eta_2$ has to be reduced with increasing problem size, which is unfavorable since it would mean to restrict the variability of the adaptive constraint set.

Both issues in theory can be tackled by restricting the original QVI to the subspace $\mathcal{N}^\perp(\mathrm{A})$ and switching from the standard Euclidean norm $\|\cdot\|_B := \|\cdot\|_2$ to the problem-specific norm $\|x\|_B := \sqrt{x^\top \mathrm{A}^\top \mathrm{A}\, x}$. We describe this approach in detail below. For practical applications this approach would require a singular value decomposition (SVD) of the operator $\mathrm{A}^\top \mathrm{A}$, which is intractable for large problem sizes.

*Restriction to $\mathcal{N}^\perp(\mathrm{A})$* . In order to deal with the missing strong monotonicity of operator $\nabla F$, we restrict the problem (3.12) to the complement $\mathcal{N}^\perp(\mathrm{A})$ of the null space $\mathcal{N}(\mathrm{A})$ of operator A:

$$\begin{aligned} &\text{Find } p^* \in \Pi_{\mathcal{N}^\perp(\mathrm{A})}(\mathcal{D}(p^*)) \text{ such that} \\ &\langle \nabla F(p^*), p - p^*\rangle \geq 0, \quad \forall p \in \Pi_{\mathcal{N}^\perp(\mathrm{A})}(\mathcal{D}(p^*)). \end{aligned} \tag{4.9}$$

This restriction is justified by the following proposition.

Proposition 4.4. *Assume that the set $\mathcal{D}(p)$ depends only on $p_{res} := \Pi_{\mathcal{N}^\perp(\mathrm{A})}(p)$.*

*(i) Let $p^*_{res}$ be a solution to the restricted problem (4.9). Then, any $p^* \in \Pi^{-1}_{\mathcal{N}^\perp}(p_{res}) \cap \mathcal{D}(p_{res})$ is a solution to the original problem (3.12).*

*(ii) If $p^*$ is a solution of the unrestricted problem, then $p^*_{res} := \Pi_{\mathcal{N}^\perp(\mathrm{A})}(p^*)$ is a solution of the restricted problem. Thus, we can express any solution $p^*$ of the unrestricted problem (3.12) as $p^* = p^*_{res} + p^*_{\mathcal{N}}$ with $p^*_{\mathcal{N}} = \Pi_{\mathcal{N}(\mathrm{A})}(p^*)$.*

The proof of Prop. 4.4 can be found in the Appendix A.

We remark that in our applications $\mathcal{D}(p)$ depends on $L^\top u = L^\top M^{-1}(f - \mathrm{A}\, p)$, thus the assumption of Prop. 4.4 is satisfied.

*Choosing a Problem Specific Norm.* We now address the issue, that the condition number $\gamma_2$ of $\nabla F$ w.r.t. the standard Euclidean norm increases with the problem size.

In order to show uniqueness of a solution to (4.9), we consider the space $\mathbb{R}^{mn} \cap \mathcal{N}^\perp(\mathrm{A})$ equipped with the norm

$$\|x\|_B := \sqrt{x^\top \mathrm{A}^\top \mathrm{A}\, x}, \tag{4.10}$$

which indeed is a norm on $\mathcal{N}^\perp(\mathrm{A})$.

Standard calculus then shows that on the subspace $\mathcal{N}^\perp(\mathrm{A})$ equipped with norm $\|\cdot\|_B$ the Lipschitz constant $\mu_B$ and monotonicity constant $\nu_B$ become 1 and uniqueness is obtained if $\eta_B < 1$.

Two open problems remain, rendering the above approach, the restriction to $\mathcal{N}^\perp(\mathrm{A})$ together with a problem specific norm, a purely academic one:

- The condition $\eta_B < 1$ is hard to verify in practice, since the projection $\Pi_{\mathcal{D}(p)}$ is defined w.r.t. the specific norm $\|\cdot\|_B$ and a closed form for this projection in general is not at hand, even if it is available for the Euclidean norm (as for example for the standard TV semi-norm).
- Restriction to the space $\mathcal{N}^\perp(\mathrm{A})$ requires the SVD of A, which, for larger images numerically is intractable.

As a consequence of these two open problems, we follow an alternative ansatz. We will see that in this ansatz the subspace $\mathcal{N}^\perp(\mathrm{A})$ and the norm $\|.\|_B$ will play also an important role.

**4.2.2. Uniqueness Results for the Proposed Approach.** We recall that the theory from Sect. 4.1 provides existence of a solution $p^*$ to the problem (3.12) under Assumption 4.1. Our considerations in the previous subsection showed that we cannot expect uniqueness of the component $p^*_{\mathcal{N}}$ of $p^*$ in the null space $\mathcal{N}$ of operator A. On the other hand, we are mainly interested in $u^* = M^{-1}(f - \mathrm{A}\, p^*)$ ($M = \mathrm{Id}$ in the case of denoising and $M$ being the blur operator in the case of deblurring), which does not depend on $p^*_{\mathcal{N}}$. We therefore focus on $v^* := \mathrm{A}\, p^*$, for which we will show uniqueness, and on the standard Euclidean norm of $v^*$. Note that the mapping $p \mapsto v := \mathrm{A}\, p$ implicitly depends only on the component of $p$ in $\mathcal{N}^\perp(\mathrm{A})$ and that $\|v\|_2 = \|p\|_B$ holds, showing the relationship to the previous section.

We prove uniqueness of $v^*$ under the following assumption:

Assumption 4.5 (for uniqueness).

*(i)* $F(p) = \frac{1}{2}\|Ap - f\|^2$. *In particular,* $\nabla F(p) = A^\top(Ap - f)$ *is Lipschitz-continuous with Lipschitz-constant* $\mu_2 := \|A^\top A\|_2$.

*(ii) The set* $\mathcal{D}(p)$ *depends only on* $v = \mathrm{A}\,p$, *i.e. there exists* $\tilde{\mathcal{D}}(\cdot)$ *such that* $\mathcal{D}(p) = \tilde{\mathcal{D}}(v) = \tilde{\mathcal{D}}(\mathrm{A}\,p)$. *In particular, instead of the variation rate of* $p \rightrightarrows \mathcal{D}(p)$, *we consider the variation rate* $\tilde{\eta}$ *of* $v \to \tilde{\mathcal{D}}(v)$, *i.e. the smallest number, such that*

$$\|\Pi_{\tilde{\mathcal{D}}(v)} q - \Pi_{\tilde{\mathcal{D}}(v')} q\|_2 \leq \tilde{\eta}\|v - v'\|_2, \quad \forall q \in \mathbb{R}^{mn}. \tag{4.11}$$

*(iii) The variation rate* $\tilde{\eta}$ *is less than* $\frac{1}{\sqrt{\mu_2}}$.

Before showing uniqueness, let us first define the operator $T(v) : \mathbb{R}^n \rightrightarrows \mathbb{R}^{mn}$ as follows: Let $\overline{p} \in T(v)$ if and only if $\overline{p} \in \tilde{\mathcal{D}}(v)$ and it is a solution to the VI

$$\langle \nabla F(\overline{p}), p - \overline{p}\rangle \geq 0 \quad \forall p \in \tilde{\mathcal{D}}(v). \tag{4.12}$$

We remark that due to our special choice of $F$, for $\tilde{\mathcal{D}}(v)$ being convex, closed and non-empty, the operator $\mathrm{A} \circ T$ is single-valued due to the strict convexity of $\min_{\tilde{v} \in \mathrm{A}\,\tilde{\mathcal{D}}(v)} \frac{1}{2}\|\tilde{v} - f\|_2^2$. We find that for any solution $p^*$ to QVI (3.12) $v^* = \mathrm{A}\,p^*$ is a fixed-point of $\mathrm{A} \circ T$.

Theorem 4.6 (Uniqueness).

1. *Under Assumptions 4.1 and 4.5 (i)-(ii), the mapping* $\mathrm{A} \circ T$ *is Lipschitz-continuous with constant* $\lambda_2 := \tilde{\eta}\sqrt{\mu_2}$.
2. *Let* $p^*$ *be a solution of the QVI* (3.12) *(cf. Prop. 4.2). If in addition Assumption 4.5 (iii) holds, then* $v^* := \mathrm{A}\,p^*$ *is unique.*

The proof of Thm. 4.6 can be found in the Appendix B.

Remark 4.7.

- *An equivalent condition to Assumption 4.5 (iii) is that the Lipschitz-constant of* $v \to \mathrm{A} \circ \Pi_{\mathcal{D}(v)}(q)$ *is less than* 1 *for all* $q$ *(cf. condition on* $\eta_B$ *in Sect. 4.2.1).*
- $\|\mathrm{A} \circ T(v_1) - \mathrm{A} \circ T(v_2)\|_2$ *is equal to* $\|T(v_1) - T(v_2)\|_B$, *where* $\|.\|_B$ *with* $B = \mathrm{A}^\top \mathrm{A}$ *is the special norm on* $\mathcal{N}^\perp(\mathrm{A})$ *considered before. Thus Theorem 4.6 provides that operator* $T$ *under the said conditions is a contraction in the norm* $\|\cdot\|_B$ *(*$B = \mathrm{A}^\top \mathrm{A}$*).*

We recall that in the considered applications for image restoration (cf. Sect. 6 and previous work [23, 24]), we are actually interested in the variable $u := M^{-1}(f - \mathrm{A}\,p)$. It follows from Theorem 4.6 that this variable is unique under Assumptions 4.1 and 4.5.

For specific examples of adaptive TV denoising, to guarantee uniqueness of the fixed-point problem, it remains provide a sufficiently small variation rate. The variation rate, on the other hand, is typically related to the regularization strength, as in Example 3.2 considered above. We revisit this example in the following:

Example 4.8. *We revisit the adaptive TV regularization in Example 3.2, where the local constraint set* $\mathcal{D}_i^{loc}$ *is given as*

$$\mathcal{D}_i^{loc} = B_{\alpha_i}(0), \qquad \alpha_i := \max\{\alpha_0(1 - \kappa(|(\mathrm{L}^\top u)_i|), \varepsilon\}. \tag{4.13}$$

*Due to the relation $u = M^{-1}(f - v)$, each $\alpha_i$ depends on $v$ by*

$$\alpha_i(v) = \max\left\{\alpha_0(1 - \kappa(|(\mathrm{L}^\top M^{-1}(f-v))_i|), \varepsilon\right\} \tag{4.14}$$
$$= \max\left\{\alpha_0(1 - \kappa(|(\mathrm{A}^\top(f-v))_i|), \varepsilon\right\}, \tag{4.15}$$

*where we used $\mathrm{A} = M^{-\top}\,\mathrm{L}$. We calculate the variation rate $\tilde{\eta}$ of $\tilde{\mathcal{D}}(v)$. Let $v, \tilde{v} \in \mathbb{R}^n, q \in \mathbb{R}^{mn}$ be arbitrary. Since the projection of $q$ onto $\tilde{\mathcal{D}}(v)$ is a scaling of the $n$ components $q_i \in \mathbb{R}^{mn}$ to at most length $\alpha_i(v)$, we find*

$$\|\Pi_{\tilde{\mathcal{D}}(v)} q - \Pi_{\tilde{\mathcal{D}}(v)} q\|_2^2 \leq \sum_{i=1}^n |\alpha_i(v) - \alpha_i(\tilde{v})|^2 \tag{4.16}$$
$$\leq \alpha_0^2\kappa^2 \sum_{i=1}^n \big| |(\mathrm{A}^\top(f-v))_i| - |(\mathrm{A}^\top(f-\tilde{v}))_i| \big|^2 \tag{4.17}$$
$$\leq \alpha_0^2\kappa^2 \sum_{i=1}^n |(\mathrm{A}^\top(f-v) - \mathrm{A}^\top(f-\tilde{v}))_i|^2 \tag{4.18}$$
$$= \alpha_0^2\kappa^2 \|\,\mathrm{A}^\top(v-\tilde{v})\|_2^2 \tag{4.19}$$
$$\leq \alpha_0^2\kappa^2 \|\,\mathrm{A}\,\|_2^2\, \|v-\tilde{v}\|_2^2. \tag{4.20}$$

*Thus $\tilde{\eta} = \alpha_0\kappa\|\,\mathrm{A}\,\|_2 = \alpha_0\kappa\sqrt{\mu_2}$. Theorem 4.6 therefore guarantees a unique solution if*

$$\alpha_0\kappa\mu_2 < 1. \tag{4.21}$$

*Considering the task of* denoising, *where $M = \mathrm{Id}$, $\mathrm{A} = \mathrm{L}$ and $\|L\|_2^2 = \mu_2 = 8$, condition (4.21) becomes $\alpha_0\kappa < \frac{1}{8}$. Given a fixed maximal regularization strength $\alpha_0$ we thus can determine feasible values for $\kappa$ to guarantee uniqueness of the solution.*

*For the task of* deblurring, *where $\mathrm{A} = M^{-\top}\,\mathrm{L}$, we expect that in practical applications $\mu_2 = \|M^{-\top}L\|_2^2 \gg 1$ due to small eigenvalues of $M$ and thus that uniqueness can be guaranteed only for very small $\alpha_0$ (weak smoothing) or $\kappa$ (weak adaptivity).* ♦

## 5. Numerics.

Throughout this section, we assume that Assumptions 4.1 and 4.5 are satisfied.

In particular, we assume that the dependency of $\mathcal{D}(p)$ on $p$ is actually a dependency on $v := \mathrm{A}\,p$. We change the notation accordingly by writing $\mathcal{D}(v)$ instead of $\mathcal{D}(p)$.

### 5.1. Proposed Algorithm.

In the following, we propose an algorithm to solve the QVI (3.12). This algorithm builds on the ideas already presented in [24]. However, we now provide convergence results for the more general case $\tilde{\eta}\gamma_2 < 1$ and, in particular, for arbitrary image sizes.

As already proposed in [24, 27], we consider an outer and an inner loop. In the outer loop we update the value $v$ which defines the constraint set $\mathcal{D}(v)$. The inner step consists in solving the variational inequality

$$\begin{aligned} &\text{Find } \overline{p} \in \mathcal{D}(v) \text{ such that} \\ &\langle \nabla F(\overline{p}), p - \overline{p}\rangle \geq 0, \quad \forall p \in \mathcal{D}(v) \end{aligned} \tag{5.1}$$

with fixed constraint set $\mathcal{D}(v)$. Recall that the operator which maps $v$ to an *exact* solution $\overline{p}$ of (5.1) is denoted by $T(v)$.

Several methods have been proposed to numerically solve (5.1). At this point, we consider some arbitrary method and denote its *numerical* result by $sol(\mathcal{D}, p^0, N)$, where $\mathcal{D}$ is the current constraint set, $p^0$ is an initial value and $N$ is the number of inner iteration steps. We assume that an a-priori error bound for this method is available: for any $\varepsilon > 0$ we can find $N$ large enough and independent of $p^0$ and $\mathcal{D}$, such that the inner problem can be solved up to an error

$$\| \mathrm{A}\big(sol(\mathcal{D}(v), p^0, N) - T(v)\big)\|_2 \leq \varepsilon, \quad \forall p^0 \in B_R(0), v \in \mathbb{R}^{mn}, \tag{5.2}$$

where $R$ is the a-priori bound on $p$ (cf. Remark 4.3). Exemplary methods fulfilling these requirements are discussed in Sect. 5.1.2.

**Algorithm 1:** Outer Iteration

**Output**: Sequence $(p^{[k]})_k$ converging to a solution $p^*$ of (3.12).
Choose arbitrary $p^{[0]} \in B_R(0) \subset \mathbb{R}^{mn}$, $v^{[0]} = \mathrm{A}\, p^{[0]}$. `// initialization`
**begin**
  **for** $k = 0, \ldots, K-1$ **do**
    $p^{[k+1]} = sol(\mathcal{D}(v^{[k]}), p^{[k]}, N)$ `// →  Algorithm solving VI (5.1)`
    $v^{[k+1]} = \mathrm{A}\, p^{[k+1]}$

**5.1.1. The Outer Iteration.** In Algorithm 1 we outline the outer iteration, which provides a sequence $p^{[k]}$ converging to a fixed-point $p^*$ of (3.12). For each iterate $p^{[k]}$ we set $v^{[k]} := \mathrm{A}\, p^{[k]}$ and fix the constraint set $\mathcal{D}(v^{[k]})$. The corresponding inner problem (5.1) is solved in an inner iteration to obtain $p^{[k+1]}$.

**5.1.2. Solving the Inner Problem.** In order to solve the inner problem (5.1) or its equivalent saddle point formulation, several approaches providing the required error estimate (5.2) have been proposed in literature. Among them are, e.g., Nesterov's method in [26], FISTA [3], and the primal-dual algorithms proposed by Chambolle & Pock [9]. Out of these candidates we exemplarily pick FISTA (with constant step size), see Algorithm 2. We refer to the iteration within the FISTA algorithm as the inner iteration. In order to distinguish the inner iterates from the outer ones, i.e. $p^{[k]}$ and $v^{[k]}$, we use the notation $p^{(k)}$ with parentheses.

We briefly recall the convergence results for FISTA [3], for which an error bound of the form (5.2) is available. We remark that similar estimates hold for the primal-dual algorithms.

**Algorithm 2:** FISTA.

**Input**: $N \geq 1$, constraint set $\mathcal{D}(v^{[k]})$, initial guess $p^{(0)} \in B_R(0) \subset \mathbb{R}^{mn}$
**Output**: $p = p^{(N)}$
Let $\tau^{(0)} = 1$, $q^{(1)} = p^{(0)}$
**begin**
  **for** $l = 0, \ldots, N-1$ **do**
    $p^{(l)} := \Pi_{\mathcal{D}(v^{[k]})}\big(q^{(l)} - \frac{1}{\mu_2}(\mathrm{A}^\top \mathrm{A}\, q^{(l)} - \mathrm{A}^\top f)\big)$
    $\tau^{(l+1)} := \frac{1}{2}(1 + \sqrt{1 + 4(\tau^{(l)})^2})$
    $q^{(l+1)} := p^{(l)} + \frac{\tau^{(l)}-1}{\tau^{(l+1)}}(p^{(l)} - p^{(l-1)})$

Lemma 5.1. *For the result obtained by FISTA applied to the problem (5.1), we have the following error estimate:*

$$\|\mathrm{A}\, p^{(N)} - \mathrm{A}\, T(v^{[k]})\|_2 \leq \frac{2\sqrt{2\mu_2}}{N+1}\|p^{(0)} - T(v^{[k]})\|_2. \tag{5.3}$$

*Using the boundedness of $p^{(0)}$ and $T(v^{[k]})$ (cf. Remark 4.3), it follows that*

$$\|\mathrm{A}\, p^{(N)} - \mathrm{A}\, T(v^{[k]})\|_2 \leq \frac{4R\sqrt{2\mu_2}}{(N+1)}. \tag{5.4}$$

*Proof.* Recall that the inner problem (5.1) is equivalent to

$$\min_{p \in \mathcal{D}(v^{[k]})} F(p), \tag{5.5}$$

where $F(p) = \frac{1}{2}\|\mathrm{A}\, p - f\|_2^2$. Let $\overline{p} = T(v^{[k]})$ denote a solution to (5.5). The inequality (5.3) is obtained from the error estimate

$$F(p^{(N)}) - F(\overline{p}) \leq \frac{2\mu_2}{(N+1)^2}\|p^{(0)} - T(v^{[k]})\|_2^2, \tag{5.6}$$

cf. Beck & Teboulle's Thm. 4.4. in [3], and

$$\frac{1}{2}\|\mathrm{A}\, p^{(N)} - \mathrm{A}\,\overline{p}\|_2^2 \leq F(p^{(N)}) - F(\overline{p}), \tag{5.7}$$

cf. Lemma C.1 in the Appendix. Combining (5.6) and (5.7) then shows (5.3). The inequality (5.4) follows from (5.3) using

$$\|p^{(0)} - T(v^{[k]})\|_2 \leq \|p^{(0)}\|_2 + \|T(v^{[k]})\|_2 \leq 2R. \tag{5.8}$$

■

In view of the next subsection, we consider the following special case.

Remark 5.2. *Assume that $T(v^{[k]}) \in \mathcal{N}^\perp(\mathrm{A})$ and that starting with a value $p^{(0)} \in \mathcal{N}^\perp(\mathrm{A})$ the sequence $p^{(k)}$ stays in this subspace. Using the basic fact that $\|x\|_2 \leq \frac{1}{\sqrt{\nu_2}}\|\mathrm{A}\,x\|_2$ for $x \in \mathcal{N}^\perp(\mathrm{A})$, where $\nu_2$ is the smallest positive eigenvalue of $\mathrm{A}^\top \mathrm{A}$, it follows from (5.3), that*

$$\|\mathrm{A}\,p^{(N)} - \mathrm{A}\,T(v^{[k]})\|_2 \leq \frac{2\sqrt{2\gamma_2}}{N+1}\|\mathrm{A}\,p^{(0)} - \mathrm{A}\,T(v^{[k]})\|_2, \tag{5.9}$$

*where $\gamma_2 = \frac{\mu_2}{\nu_2}$ is the condition number of $\mathrm{A}^\top \mathrm{A}$ restricted to $\mathcal{N}^\perp(\mathrm{A})$.*

**5.2. Convergence.** In the following, we show convergence of the proposed Algorithm 1 and provide convergence rates for a special case.

Proposition 5.3 (Convergence). *Let Assumptions 4.1 and 4.5 be satisfied. Moreover, assume that $sol(\mathcal{D}(v), p, N)$ provides an approximate solution of (5.1) with an error less than $\varepsilon > 0$ (independent from $p \in B_R(0)$), i.e.*

$$\|\mathrm{A}\,sol(\mathcal{D}(v), p, N) - \mathrm{A}\,T(v)\|_2 \leq \varepsilon \qquad \forall p \in B_R(0). \tag{5.10}$$

*Then, the following holds:*

*(i) for the solution $v^{[k]}$ of Algorithm 1 we have*

$$\|v^{[k]} - v^*\|_2 \leq \varepsilon\frac{1}{1-\lambda_2} + \lambda_2^k\|v^0 - v^*\|_2, \tag{5.11}$$

*where $v^*$ is the unique fixed-point of $\mathrm{A}\circ T$ and $\lambda_2$ is the Lipschitz-constant of $\mathrm{A}\circ T$.*

*(ii) $\{v^{[k]}\}_k$ converges to the fixed-point $v^*$ for $\varepsilon \to 0$ and $k \to \infty$.*

*Proof.* We have

$$\begin{aligned}
\|v^{[k]} - v^*\|_2 &= \|\mathrm{A}\,p^{[k]} - \mathrm{A}\,p^*\|_2 = \|\mathrm{A}\,p^{[k]} - \mathrm{A}\,T(v^*)\|_2 && (5.12)\\
&\leq \|\mathrm{A}\,p^{[k]} - \mathrm{A}\,T(v^{[k-1]})\|_2 + \|\mathrm{A}\,T(v^{[k-1]}) - \mathrm{A}\,T(v^*)\|_2 && (5.13)\\
&\leq \|\mathrm{A}\,sol(v^{[k-1]}) - \mathrm{A}\,T(v^{[k-1]})\|_2 + \lambda_2\|v^{[k-1]} - v^*\|_2 && (5.14)\\
&\leq \varepsilon + (\lambda_2\varepsilon + \lambda_2^2\|_2 v^{[k-2]} - v^*\|_2) && (5.15)\\
&\leq \varepsilon(1 + \lambda_2 + \lambda_2^2 + \dots + \lambda_2^{k-1}) + \lambda_2^k\|v^{[0]} - v^*\|_2. && (5.16)
\end{aligned}$$

Using the limit of the geometric series, we deduce claim (i). Claim (ii) follows from (i) under Assumption 4.5 (iii), since then $\lambda_2 < 1$ and thus $\lambda_2^k \to 0$ for $k \to \infty$. ■

Proposition 5.4 (Convergence rates). *Let Assumptions 4.1 and 4.5 be satisfied. Moreover, assume that the inner problem (5.1) can be solved with an error bound*

$$\|\mathrm{A}\,p^{(N)} - \mathrm{A}\,T(v^{[k]})\|_2 \leq \frac{\delta}{4}\|\mathrm{A}\,p^{(0)} - \mathrm{A}\,T(v^{[k]})\|_2, \tag{5.17}$$

*where $\delta := 1 - \tilde{\eta}\sqrt{\mu_2}$ is the contraction gap of the problem (3.12). (Recall that $T(v^{[k]})$ is the exact solution of the inner problem.) Consider a solution $p^*$ of (3.12) and $v^* := \mathrm{A}\,p^*$. Then, Algorithm 1 converges according to*

$$\|v^{[K]} - v^*\|_2 \leq \frac{1}{\delta}\exp\left(-\frac{\delta}{2}K\right)\|v^{[0]} - \mathrm{A}\,T(v^{[0]})\|_2, \tag{5.18}$$

*where $K$ is the number of outer iterations.* *Proof.* see Appendix D.

Remark 5.5.

**1D case:** *In the one-dimensional case the null space $\mathcal{N}(\mathrm{A})$ of $\mathrm{A}$ is spanned by the vector $(0,\dots,0,1)^\top$. When considering Example 4.8 and Algorithm 2 (FISTA), we can guarantee that the sequence $\{p^{(k)}\}_k$ stays in $\mathcal{N}(\mathrm{A})^\perp$, provided that the initial value $p^{(0)}$ is chosen in $\mathcal{N}(\mathrm{A})^\perp$. This is due to the fact that the projection onto $\mathcal{D}(v^{[k]})$ decouples into independent projections to 1D intervals for each coordinate. In particular, if $p_n^{(0)} = 0$, the constraint $p_n^{(0)} \in [-\alpha_n, \alpha_n]$ is fulfilled and thus $p_n^{(0)}$ is not changed by the projection $\Pi_{\mathcal{D}(v^{[k]})}(p)$, i.e. $p_n^{(0)} = p_n^{(1)} = \dots = p_n^{(N)} = 0$. The sequence $p^{(k)}$ then converges to the solution $\overline{p} \in T(v^{[k]}) \cap \mathcal{N}(\mathrm{A})^\perp$, which is unique in this subspace. As a consequence, starting with $p^{[0]} \in \mathcal{N}(\mathrm{A})^\perp$ for the* outer *iteration, we can guarantee also $p^{[k]} \in \mathcal{N}(\mathrm{A})^\perp$, such that the initial value of every subsequent inner problem is again in $\mathcal{N}^\perp(\mathrm{A})$. The error estimates for FISTA (cf. Remark (5.2)) then provide the necessary conditions to apply Prop. 5.4.*

**2D case:** *In the two-dimensional case there exist counter examples indicating that (5.17) does not hold in general. The reason is that the convergence depends on the component $\|p^{(0)} - \overline{p}\|_{\mathcal{N}(\mathrm{A})}$ of the initial error $\|p^{(0)} - \overline{p}\|_2$ ($\overline{p} \in T(v^{[k]})$ fixed).*

## 6. Experiments.

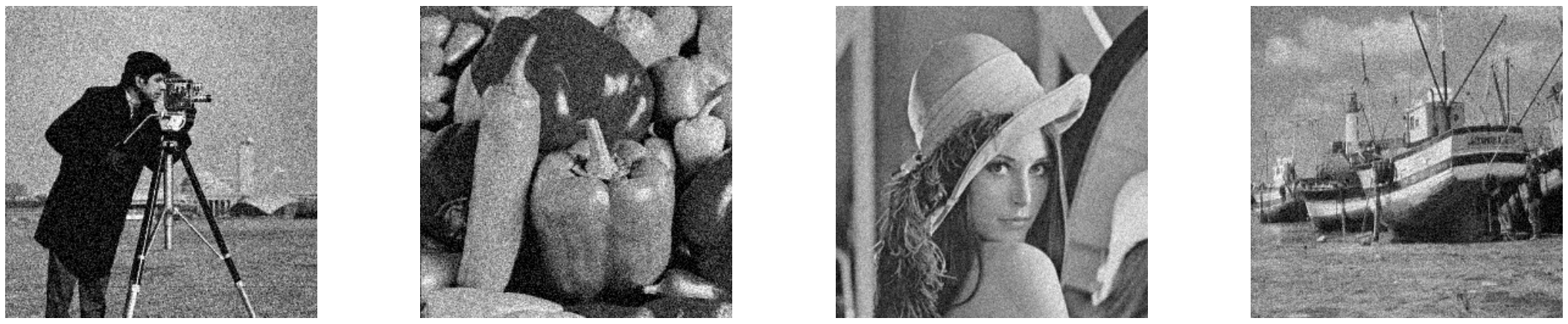

(a) Cameraman (b) Peppers (c) Lena (d) Boat

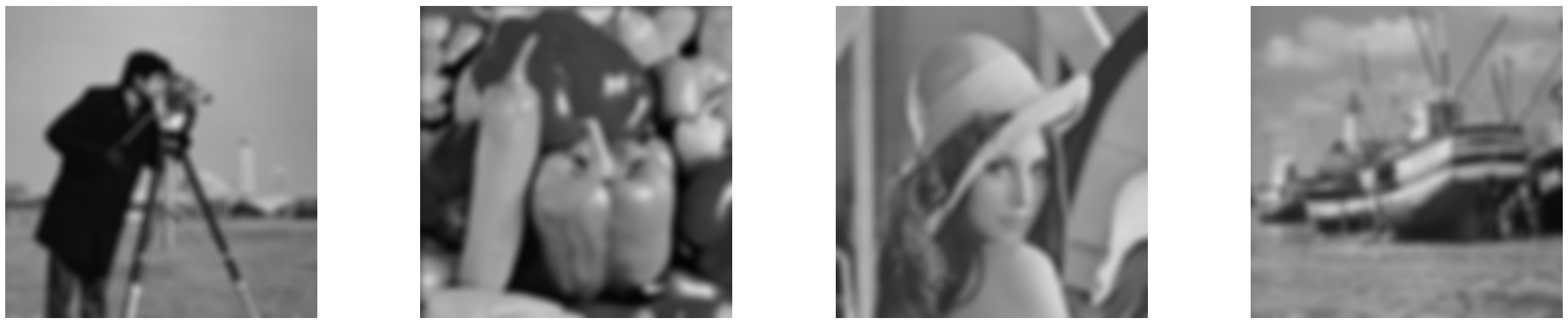

(e) Cameraman (f) Peppers (g) Lena (h) Boat

**Figure 1.** *Test images used for evaluation. Top row: images with Gaussian noise (zero mean, standard deviation* 0.1*) to evaluate denoising. Bottom row: blurred images with Gaussian noise (zero mean, standard deviation* 0.01*) to evaluate deblurring.*

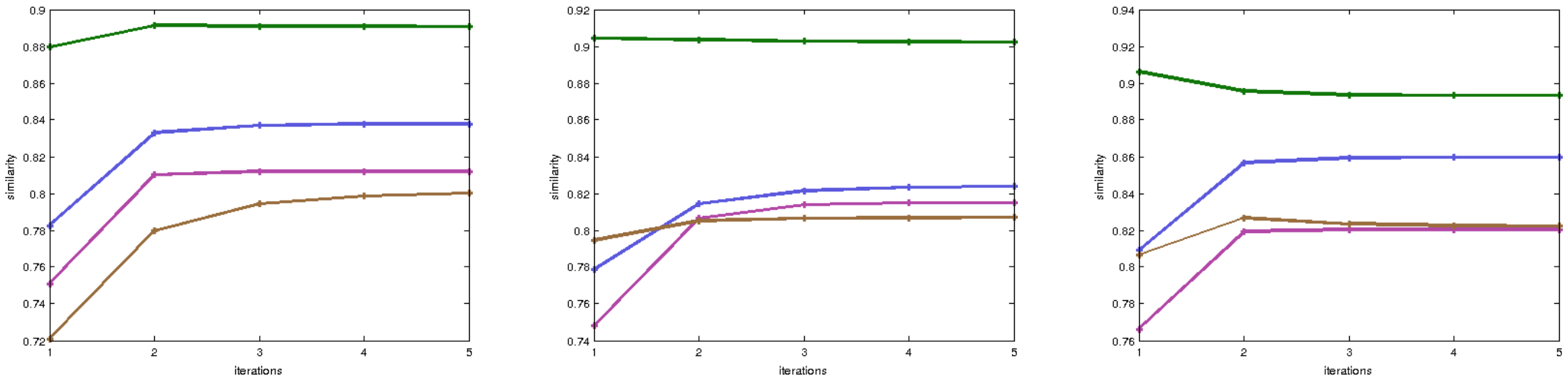

(a) Adaptive first-order (b) Anisotropic first-order (c) Adaptive 1st & 2nd order

**Figure 2.** Image denoising. (Best viewed in color.) *Evolution of the similarity measure* MSSIM *[39] between the current iterate and the ground truth image during the outer iteration of our approach with three different regularization terms. Test images: cameraman (purple), peppers (blue), Lena (green) and boat (brown). In most cases, the similarity increases during the outer iteration. This shows that the adaptivity improves by switching from a data-driven to a solution-driven model.*

**6.1. Improvement by Solution-Driven Adaptivity.** In the following, we demonstrate the benefits of applying our solution-driven adaptivity compared to the data-driven variant. To this end, we consider four different standard test images, the *cameraman*, *peppers*, *Lena* and the *boat* image, which are scaled to the range $[0, 1]$. From each image, we generate test data for the denoising problem by adding Gaussian noise with zero mean and standard deviation 0.1, and for the deblurring problem by applying a blurring operation and adding Gaussian noise with zero mean and standard deviation 0.01. The resulting images are show in Fig. 1. On these test images we evaluate the three different adaptive regularizations presented in Sect. 2.3, namely adaptive first-order TV regularization (Example 2.1), anisotropic first-order TV regularization (Example 2.2) and adaptive first- and second-order regularization (Example 2.3).

For denoising, we use the input data to initialize the constraint set $\mathcal{D}(v^{[0]})$, $v^{[0]} := f$. Therefore, running the algorithm with only one outer iteration implements a data-driven adaptivity, while running it with more than one outer iteration gives a solution-driven adaptivity. We set the required parameters to obtain a suitable result for the data driven approach and apply the solution driven variants with five outer iterations. To quantitatively evaluate the results, we make use of the similarity measure $\mathrm{MSSIM}(a, b)$ for two images $a$ and $b$ proposed by Wang et al. [39]. This measure is well suited in particular to compare restored images with their ground truth, since it is sensitive to remaining distortions.

The evolution of $\mathrm{MSSIM}(u_{sd}^{[k]}, u_{orig})$, where $k$ is the index for the outer iteration, for the three kinds of adaptive regularization and each test image is depicted in Fig. 2. Except for two cases, the solution-adaptive regularization improves the similarity from the first to the second outer iteration (recall that $k = 1$ provides the data-driven result). In most cases the similarity stays constant or is even further improved in the subsequent iteration steps. To also give a visual impression of this improvement, we depict the respective results for the cameraman image in Fig. 3. Since the differences are best

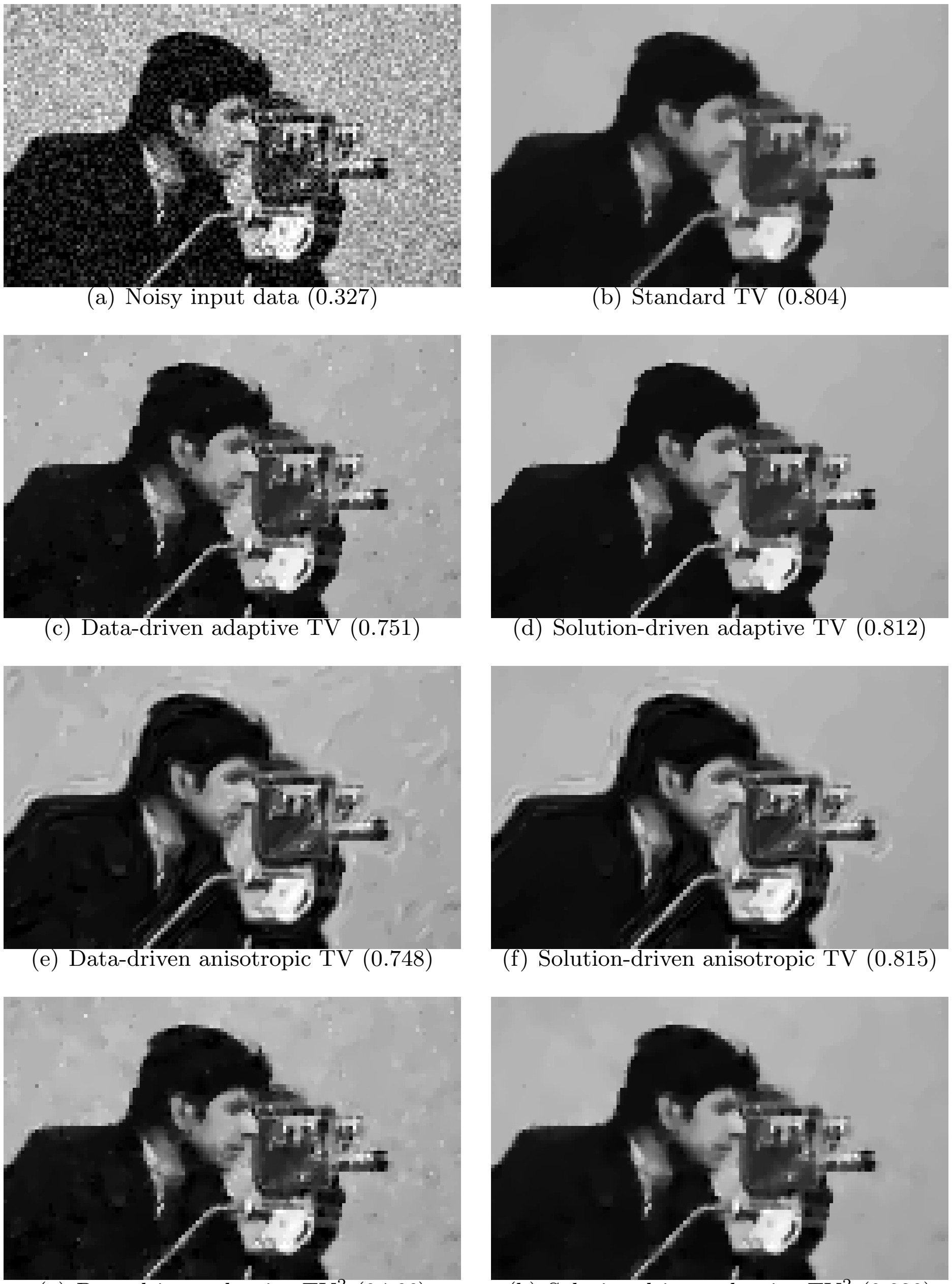

(a) Noisy input data (0.327) (b) Standard TV (0.804)

(c) Data-driven adaptive TV (0.751) (d) Solution-driven adaptive TV (0.812)

(e) Data-driven anisotropic TV (0.748) (f) Solution-driven anisotropic TV (0.815)

(g) Data-driven adaptive $\mathrm{TV}^2$ (0.766) (h) Solution-driven adaptive $\mathrm{TV}^2$ (0.820)

**Figure 3.** *Close-up of the results of* denoising *the cameraman image with different regularization approaches. The similarity values (cf. [39]) to the ground truth are given in parentheses. They correspond to those plotted in Fig. 2. The solution-driven approaches enhance the reconstruction compared to the data-driven ones and standard TV regularization. In particular, artifacts from noise are reduced.*

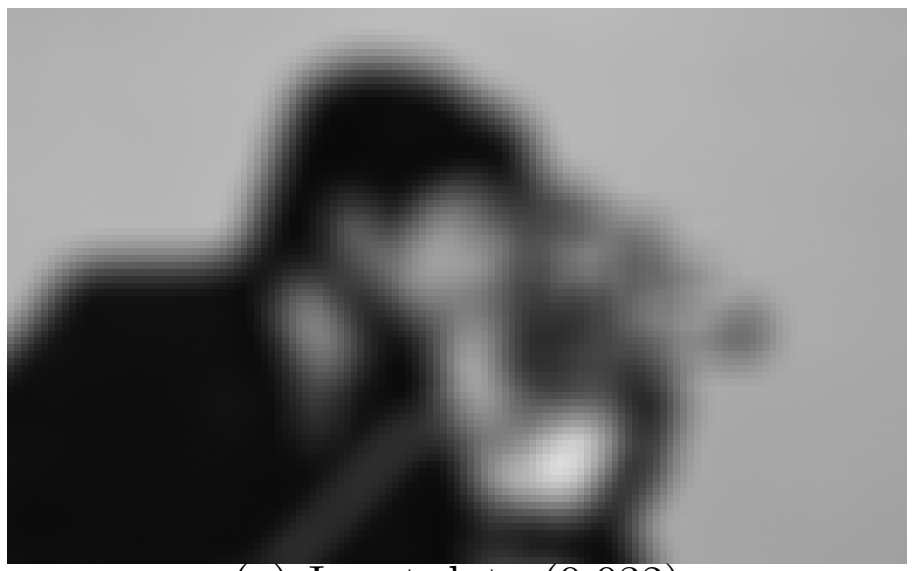

(a) Input data (0.032)

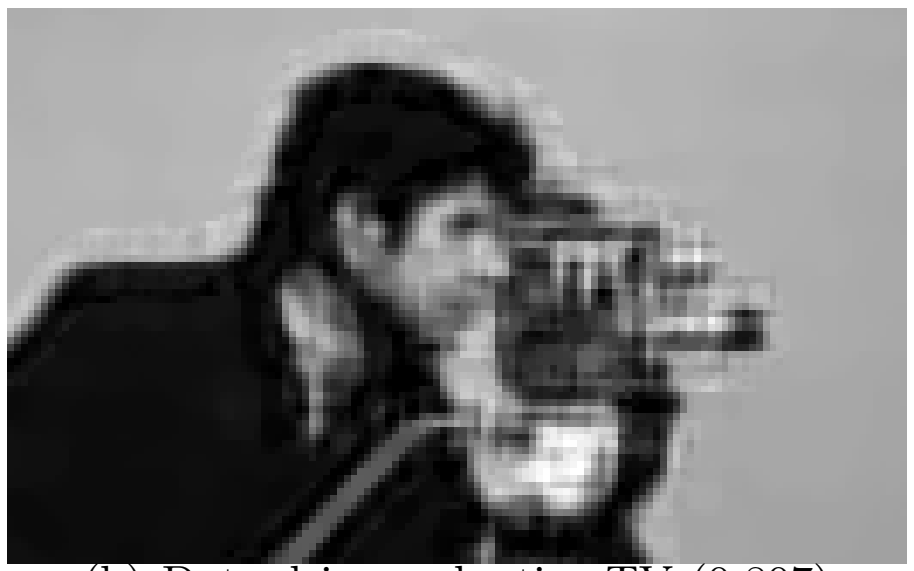

(b) Data-driven adaptive TV (0.807)

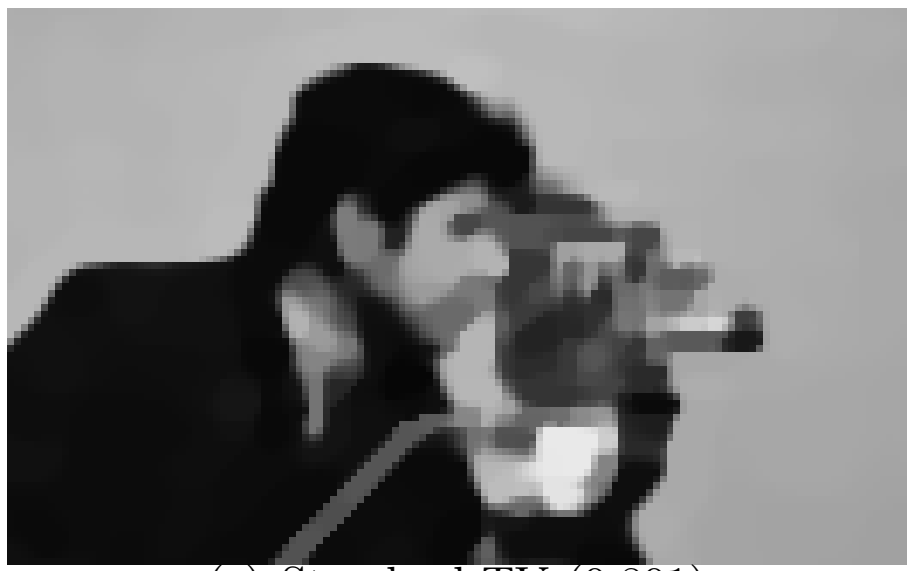

(c) Standard TV (0.801)

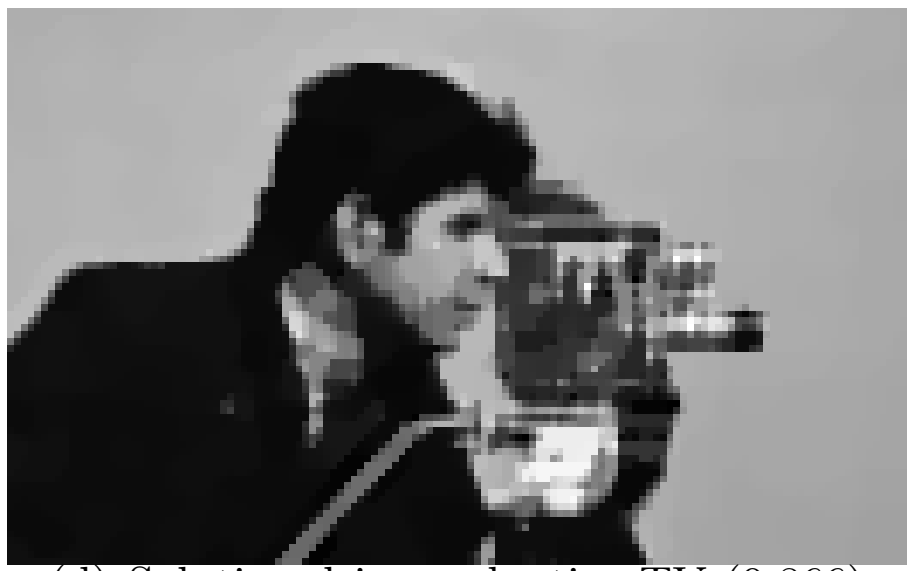

(d) Solution-driven adaptive TV (0.866)

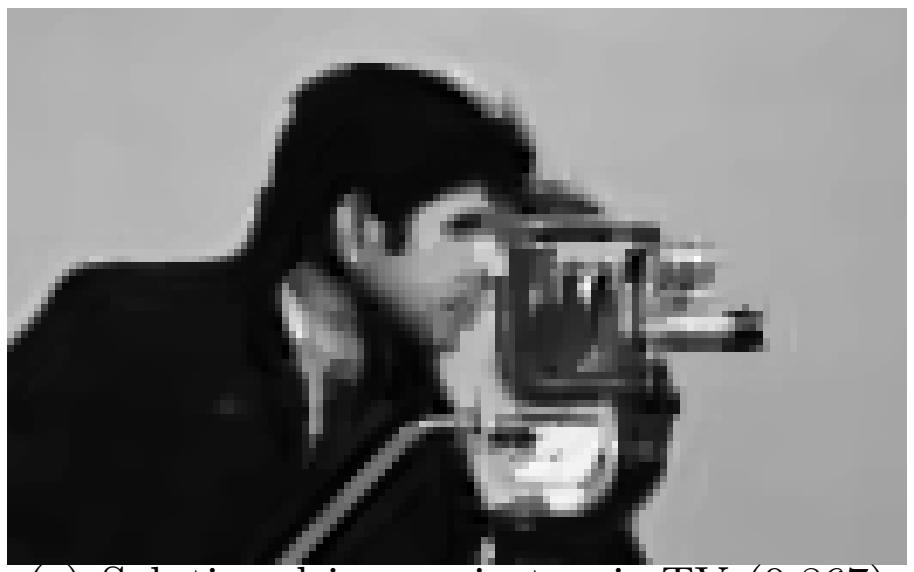

(e) Solution-driven anisotropic TV (0.867)

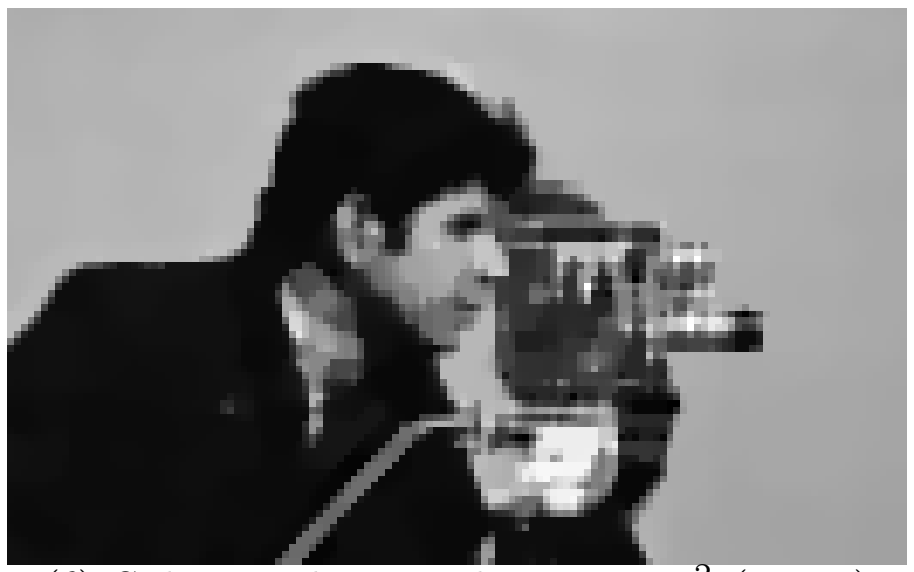

(f) Solution-driven adaptive $\mathrm{TV}^2$ (0.864)

**Figure 4.** *Close-up of the results of* deblurring *the cameraman image with different regularization approaches. Given in parentheses is the similarity to the ground truth. The data-driven approaches suffer from artifacts (e.g., in the adaptive first-order* TV *case (b)). We therefore compare our methods to standard TV (c). The solution-driven approaches enhance the reconstruction in terms of similarity to the original data compared to standard TV. Anisotropic TV regularization gives the best result.*

visible in full resolution, we focus on a close-up of the head region of the cameraman. The improvement of the similarity after five outer iterations compared to the similarity of the data-driven results are shown in Table 1 with the values in percent and averaged over the four test images. We found that using the peak-signal-to-noise-ratio (PSNR) instead of MSSIM shows a similar trend.

The theory presented in Sect. 4.2 allows us to check for each method, if the obtained result is unique and in particular independent from the initialization of the algorithm.

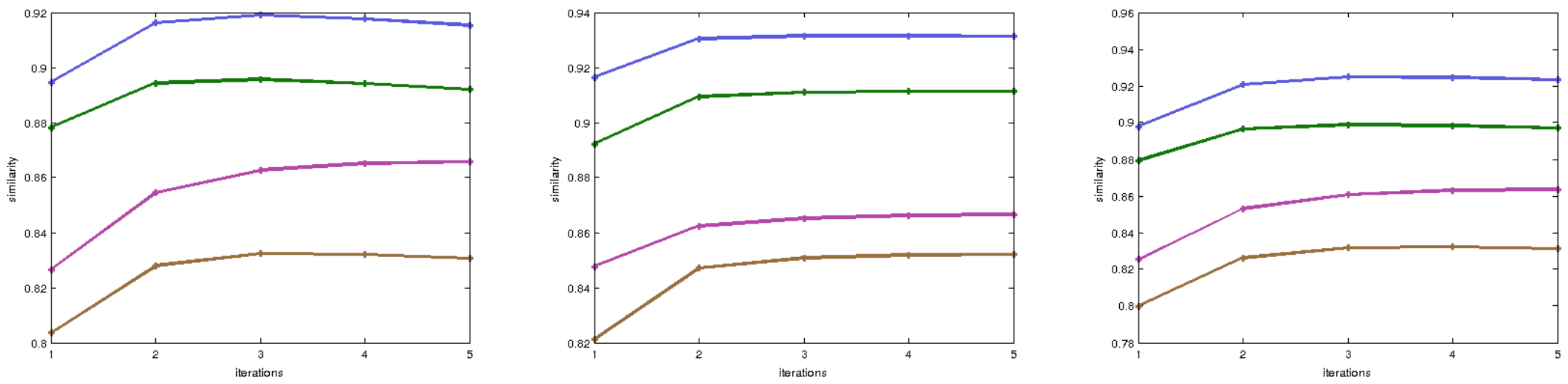

(a) Adaptive first-order (b) Anisotropic first-order (c) Adaptive 1st & 2nd order

**Figure 5.** Image deblurring.(Best viewed in color.) *Evolution of the similarity to the ground truth during the outer iteration of our approach with three different regularization terms. Test images: cameraman (purple), peppers (blue), Lena (green) and boat (brown). We observe a strong increase in similarity between the first two iteration steps. After these two steps the similarity stays almost constant, with a slight decrease in some cases. Inspecting the results visually, we observe that the results still get sharper during the later steps, which probably leads to an over-sharpening compared to the original images.*

| problem | method | gain in similarity | guaranteed uniqueness |
|---|---|---|---|
| denoising | adaptive first-order | 6.8 % | yes |
| | anisotropic first-order | 4.0 % | ? |
| | adaptive 1st & 2nd order | 3.5 % | yes |
| deblurring | adaptive first-order | 7.2 % | no |
| | anisotropic first-order | 8.8 % | no |
| | adaptive 1st & 2nd order | 4.8 % | no |

**Table 1**

*Gain in similarity to the ground truth by introducing solution-driven adaptivity for three different regularizations. We compare to the data-driven variants for denoising and to standard TV for deblurring. The values are averaged over the four test images. The theory presented in Sect. 4 allows to check for each case, if uniqueness of the result can be guaranteed, see the right column.*

Uniqueness is guaranteed for the adaptive TV regularization for first- and second-order, where the parameters $\alpha_0$ and $\kappa$ where chosen small enough to assert $\alpha_0 \kappa \mu_2 < 1$ (cf. Example 4.8). In the case of anisotropic TV regularization, it can be shown that the projection $\Pi_{\tilde{\mathcal{D}}(v)}(q)$ of $q$ onto the constraint set $\tilde{\mathcal{D}}(v)$ is Lipschitz-continuous w.r.t. $v$ (we refer to [24] for details). However, an analytic estimate of the Lipschitz constant $\tilde{\eta}$ is not at hand. Experiments indicate that $\tilde{\eta} < 0.06$ for our particular parameter setting. Assuming that this estimate is correct, our theoretical results therefore guarantee uniqueness, since $\tilde{\eta} < \frac{1}{8}$.

In the case of deblurring, it turns out that applying the data-driven approaches does not provide satisfactory results since spurious structures occur independent from method and input image (see e.g. Fig. 4(b)). Similarly, applying the solution-driven approaches using the input data as initialization results in the same artifacts. However, solution-driven approaches which start with a *constant* image as initialization provide satisfactory results (see Fig. 4(d)-(e)). This already indicates that uniqueness of the underlying QVIP can not be expected in the case of deblurring. We further comment on this below.

Evaluating with the similarity measure, see Fig. 5, we observe a substantial increase of the similarity to the ground truth during the first two outer iterations of our approach. After the second outer iteration only a slight further improvement or, in rare cases, a decrease occurs. The same effect can be also observed with the PSNR. Inspecting the results visually shows that the results actually do not become worse in terms of visual appearance, but that the image sharpness further increases, which we interpret as a slight over-sharpening of the result.

Since the data-driven approaches do not provide accurate results, we refrain from using them to compare with the solution-driven approaches. Instead, we compare to standard TV regularization with the same regularization parameter. The average improvement by solution-driven adaptivity in percent are shown in Table 1.

Concerning uniqueness of the results in the case of deblurring, we remark that, since the smallest eigenvalue of operator $M$ becomes very small (cf. Example 4.8), uniqueness can not be guaranteed for the values of $\alpha_0, \beta_0$, and $\kappa$ used in our experiments.

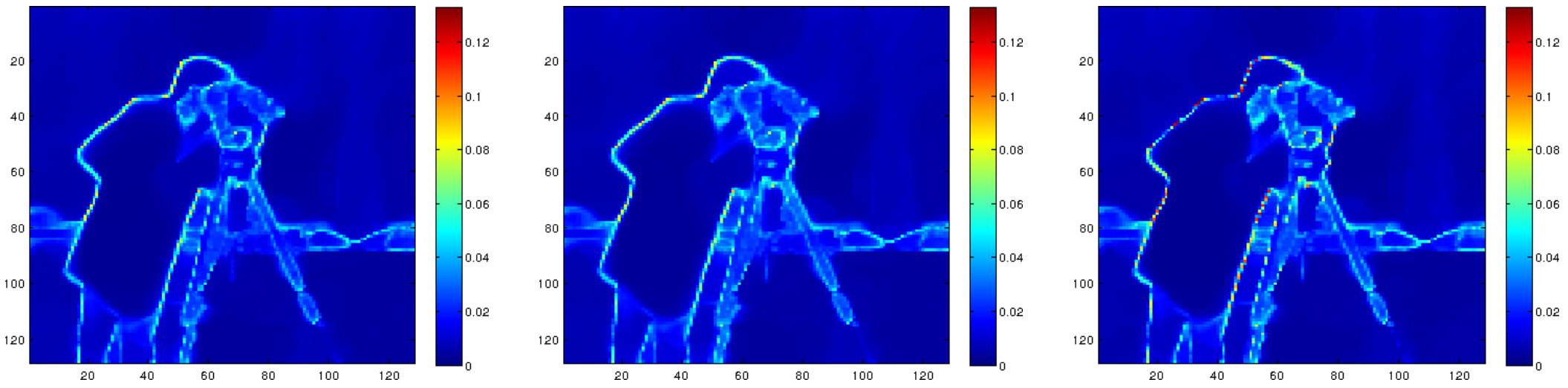

(a) Standard TV (b) Data-driven (c) Solution-driven

**Figure 6.** Dependence on the input data $f$. (Best viewed in color.) *We consider denoising of the cameraman image with the adaptive first-order TV regularization, where we denoise* 100 *versions with different realizations of additive Gaussian noise with zero mean and standard deviation* 0.1 *and calculate the pixel-wise standard deviation of the results. The maximal values attained are* 0.095 *for the standard TV approach,* 0.092 *for the data-driven and* 0.133 *for the solution-driven approach. We conclude that the sensitivity of the proposed solution-driven approach against variations in the input data is slightly higher, but of the same order of magnitude compared to non-adaptive (standard TV) and data-driven approaches.*

**6.2. Dependence on Input, Initialization and Parameters.** In the following, we discuss the dependence of the proposed algorithm on input, initialization and parameters. We focus on the adaptive TV regularization proposed in Example 2.1 in the context of image denoising.

In order to experimentally evaluate the dependence of our solution-driven approach on the input data $f$, we sample 100 noisy variants of the cameraman image with additive Gaussian noise with zero mean and standard deviation 0.1. After denoising each image, we determine the *pixel-wise* standard deviation over all 100 output images and compare our method to the data-driven variant and to standard TV. It turns out that high standard deviations for each method occur mainly along dominant edges in the image, e.g. along the silhouette of the cameraman, see Fig. 6. The maximal standard deviation attained is 0.095 for the standard TV approach, 0.092 for the data-driven and 0.133 for the solution-driven approach. We conclude for this example, that the sensitivity of our approach to variations of the input data is slightly higher, but in the same range as for the other methods.

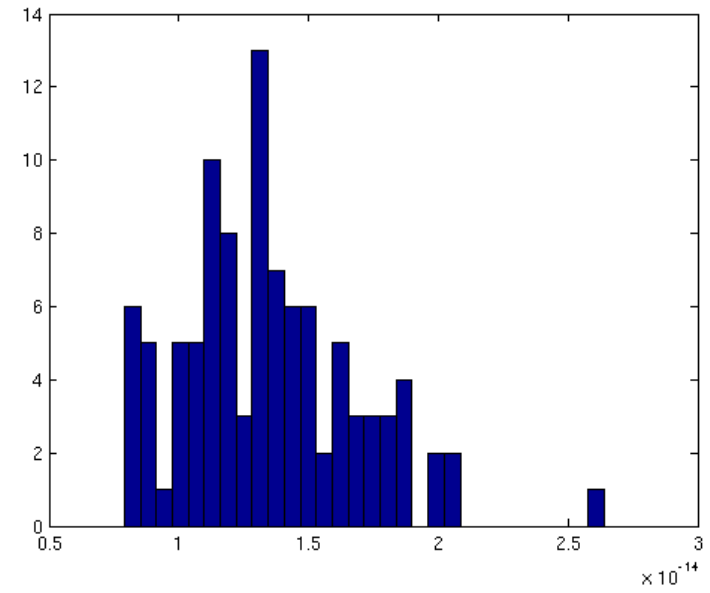

**Figure 7.** *Error histogram of 100 runs with random initialization, showing the errors after 7 outer iteration steps. The errors cannot be decreased further by increasing the number of iterations. The errors are smaller than* $3\cdot10^{-14}$ *and are fairly independent from the initialization.*

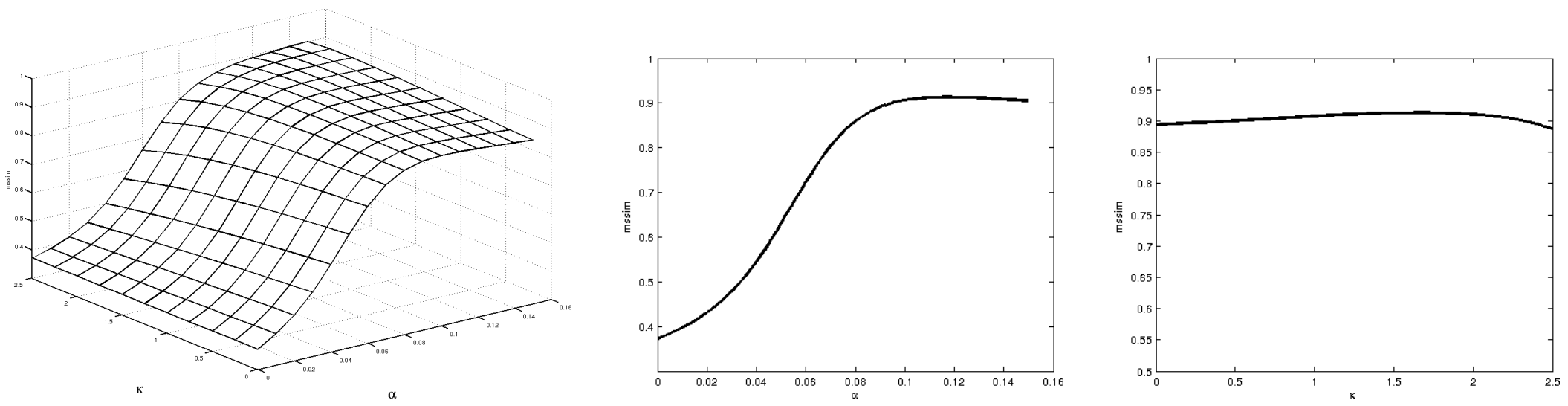

(a) (b) (c)

**Figure 8.** *Dependence on the parameters* $\alpha_0$ *and* $\kappa$ *of the similarity to the ground truth image when denoising the cameraman image with adaptive TV regularization. (a) visualization as surface* $\mathrm{MSSIM}(\alpha_0, \beta)$*. (b) cross-section through maximum along* $\alpha_0$*-axis. (c) cross-section through maximum along* $\kappa$*-axis. We observe a smooth dependence. One flat maximum (optimal parameter settings) occurs.*

Concerning initialization, our theoretical findings guarantee uniqueness of the result as long as $\alpha_0\kappa < \frac{1}{\mu_2}$, where $\mu_2 = 8$ in the case that $\mathrm{A} = \mathrm{L}$ is a discrete divergence operator. In particular, the numerical solution is independent from the initialization of the constraint set $\mathcal{D}(v^{[0]})$. However, we check this experimentally. Fig. 7 shows the distribution of the numerical error after 7 outer iterations of the proposed algorithm for 100 randomly chosen initializations $v^{[0]}$. It turns out that the error to the analytic solution is in the range of $10^{-14}$ and cannot be further decreased by additional outer iterations. Moreover, it is fairly independent from the initialization. This supports our theoretical result on uniqueness of the fixed-point.

A third issue is the dependence on the parameters $\alpha_0$ and $\kappa$. To evaluate this dependence, we run our denoising algorithm on the cameraman test image with different parameter settings $(\alpha_0, \kappa)$ taken from a grid $\{0, 0.002, \ldots, 0.15\} \times \{0, 0.05, \ldots, 2.5\}$. Since the algorithm has to be run for a large number of times, to reduce the computational effort, we restrict this experiment to the head region of the cameraman image. For each result, we evaluate its similarity to the ground truth. The resulting 2D surface $\mathrm{MSSIM}(\alpha_0, \kappa)$ is depicted in Fig. 8(a). One relatively flat maximum occurs at $(\alpha_0, \kappa) = (0.12, 1.7)$. Figs. 8(b) and (c) show cross-sections through this maximum along the $\alpha_0$- and $\kappa$-axis, respectively. Unfortunately, in this case, the optimal parameters lie outside the region where uniqueness is guaranteed. Fixing $\alpha_0 = 0.12$, the parameter $\kappa$ would need to be less than 1.08 to assert uniqueness.

From the results shown in Fig. 8 we conclude a smooth dependence on both parameters. Moreover, the flatness of the maximum guarantees robustness w.r.t. parameter variations. In practice, choosing parameters in a relatively broad neighborhood to the unknown optimal values already provides satisfactory results.

**6.3. Relation to Non-Convex Regularization.** Introducing adaptivity in TV regularization locally changes the way how the gradient (or higher derivatives) of the final solution is penalized. To gain insight into this effect, with a given solution $u^*$, one can study the empirical distribution of $|\nabla u^*(x_i)|$ versus $\langle \nabla u^*(x_i), p(x_i) \rangle$ (borrowing the notation from the continuous setting). We do this exemplarily for the case of denoising the cameraman image and adaptive TV regularization, where $\langle \nabla u^*(x_i), p(x_i) \rangle = \alpha_i \|\nabla u^*(x_i)\|$.

Studying the distribution of the norm of the discrete gradients, $|\nabla u^*(x_i)|$, versus their penalization in the regularization term, $\alpha_i |\nabla u^*(x_i)|)_i$, see Fig. 9, one recognizes that our fixed-point based approach mimics a non-convex regularizer. For the other three test images, we observe similar distributions.

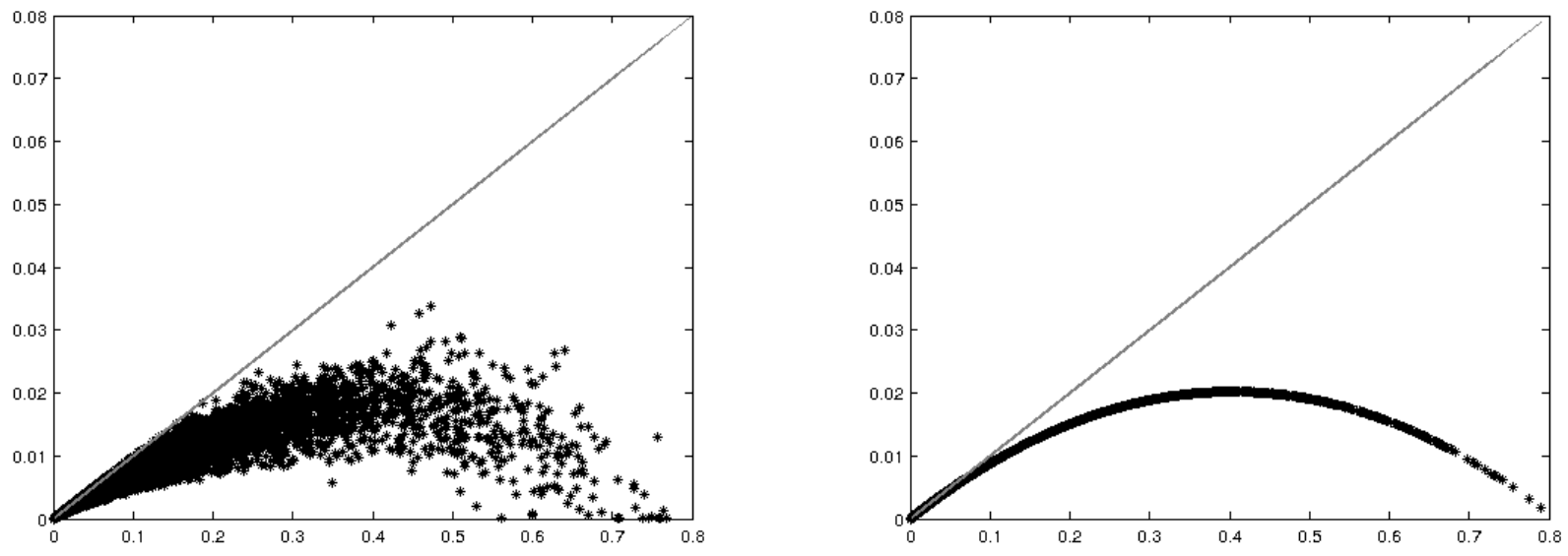

(a) Data-driven (b) Solution-driven

**Figure 9.** *A-posteriori distribution of* $(|\nabla u^*(x_i)|, \alpha_i \cdot |\nabla u^*(x_i)|)_i$ *(black dots) of (a) the data-driven* ($\alpha = \alpha(f)$) *and (b) the solution-driven result* ($\alpha = \alpha(u^*)$) *for the cameraman image, cf. Fig. 3(c) and (d), respectively. We compare these distributions to* $(|\nabla u^*(x_i)|, \alpha_0 |\nabla u^*(x_i)|)$ *(gray line), where* $\alpha_0 = 0.1$ *is the maximum value of the regularization strength. We observe that the solution driven approach mimics a non-convex regularizer.*

**6.4. Convergence.** In order to verify the convergence of our algorithm, we consider the example of an one-dimensional adaptive TV regularization, where an analytic solution can be provided. We remark that there is only a limited number of examples, for which an analytic solution for the ROF model is available. In such cases the problem (2.23) can be reformulated to a fixed-point problem in $\alpha$.

Example 6.1. *We study a discrete variant of the continuous functional in Example 3.2 with one-dimensional data. Consider a equi-distant grid of $n$ grid points. W.l.o.g. we assume the grid size to be 1. For $u, f \in \mathbb{R}^n$ let*

$$E(u) := \frac{1}{2}\sum_{i=1}^{n} |u_i - f_i|^2 + \sum_{i=1}^{n-1} \alpha_i |u_{i+1} - u_i|, \tag{6.1}$$

*where we define $\alpha \in \mathbb{R}^{n-1}$ by*

$$\alpha_i = \max(\alpha_0(1 - \kappa|u^0_{i+1} - u^0_i|), \varepsilon) \tag{6.2}$$

*for a fixed $u^0 \in \mathbb{R}^n$. Recall that we are searching for a fixed-point of $u^0 \to \arg\min_u E(u)$.*

*We consider data $f$ to be given as follows: We assume $n = 3N$ for some $N > 0$, such that the grid nodes can be divided into three disjoint sets $I_1 := \{1, \dots, N\}$, $I_2 := \{N+1, \dots, 2N\}$ and $I_3 := \{2N+1, \dots, 3N\}$. Now let*

$$f_i = \begin{cases} 0 & \text{if } i \in I_1 \cup I_3, \\ 1 & \text{if } i \in I_2. \end{cases} \tag{6.3}$$

*It can be shown that any solution of the inner problem asserts $u_i \in [0, 1]$. We make the ansatz*

$$u_i = \begin{cases} a & \text{if } i \in I_1 \cup I_3, \\ b & \text{if } i \in I_2. \end{cases} \tag{6.4}$$

*for $0 \le a \le b \le 1$. We show below that a fixed-point of this form exists. Assuming this form of $u$ and analogously for $u^0$ in (6.1), the objective function simplifies to*

$$E(u) = E(a, b) = Na^2 + \frac{N}{2}(b-1)^2 + 2\tilde{\alpha}(a^0, b^0)(b - a), \tag{6.5}$$

*where*

$$\tilde{\alpha}(a^0, b^0) := \max\{\alpha_0(1 - \kappa(b^0 - a^0)), \varepsilon\}. \tag{6.6}$$

*Standard calculus (see Appendix E) then shows, that, as long as $\kappa \le 1 - \frac{\varepsilon}{\alpha_0}$ and $\kappa\alpha_0 < \frac{N}{3}$, a fixed-point $u^*$ of $u^0 \to \arg\min_u E(u)$ of the form (6.4) is given by*

$$a^* := \frac{\tilde{\alpha}}{N}, \quad b^* := 1 - \frac{2\tilde{\alpha}}{N}, \tag{6.7}$$

*where $\tilde{\alpha} = \tilde{\alpha}(a^*, b^*) = \frac{\alpha_0(1-\kappa)}{1 - \frac{3\alpha_0\kappa}{N}}$.*

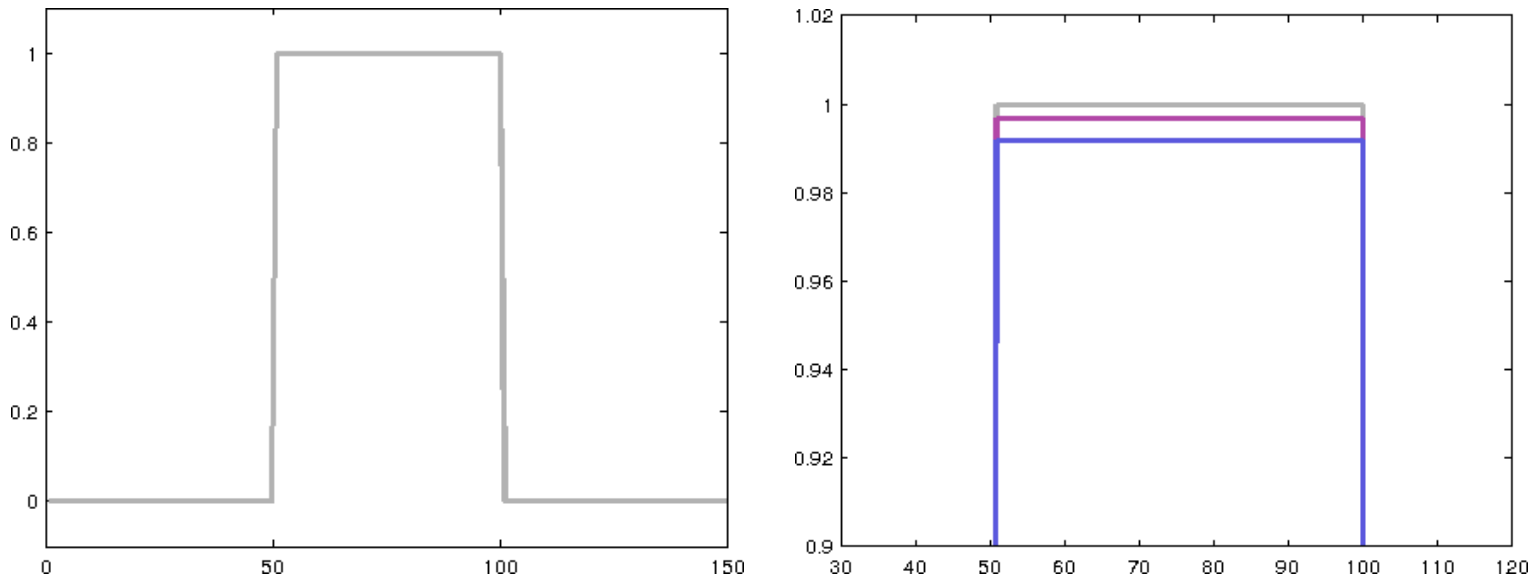

(a) data

(b) ROF and solution-driven adaptive solution

**Figure 10.** 1D example with analytic solution (cf. Example 6.1). (Best viewed in color.) *Left: data $f$. Right: solution $u^*$ of* (6.1) *(zoom, red line) for* 150 *grid nodes and parameters* $\alpha_0 = 0.2$, $\kappa = 0.6$ *compared to the solution of the standard ROF model (blue line) with* $\alpha = 0.2$. *Smoothing with the proposed approach provides a result with higher contrast compared to the standard ROF model.*

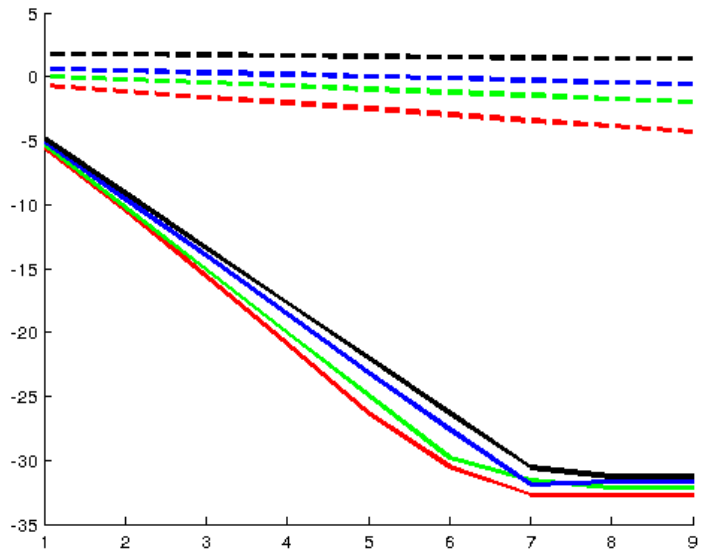

**Figure 11.** Convergence rates. (Best viewed in color.) *Plot of logarithmic numerical error (solid lines) to analytic solution over outer iteration steps for different contraction gaps* $\delta = 1 - \kappa\alpha_0\mu_2 \approx 0.1, 0.3, 0.5, 0.9$ *given by* $\alpha_0 = 2$, $\kappa = \frac{t}{4\alpha_0}$ *and* $t \in \{0.9, 0.7, 0.5, 0.1\}$ *(black, blue, green, red). The plot shows an exponential error decay, which stays well below the theoretical bound (dashed lines). The bending between step 6 and 7 is caused by the fact that the point-wise errors reach machine accuracy.*

*Note that in the one-dimensional case four is a tight upper bound for* $\mu_2$. *Thus, the condition* $\kappa\alpha_0 < \frac{1}{4}$ *to assert Assumption 4.5 (iii) independent from the problem size is sufficient to guarantee* $\kappa\alpha_0 < \frac{N}{3}$. *Moreover, our theoretical findings show that* $u^*$ *given by* (6.7) *is the* unique *fixed-point.* ♦

By means of Example 6.1, we experimentally verify the convergence rate provided by Prop. 5.4. To this end, we solve the corresponding QVI numerically with the proposed algorithm. Fig. 11 shows the theoretical and experimental convergence rates (logarithmic error over time steps) for this example and different contraction gaps $\delta = 1 - \lambda_2 = 1 - \alpha_0\kappa\mu_2$. The experimental errors $\|u^{[K]} - u^*\|_2 = \|\operatorname{L} p^{[K]} - \operatorname{L} p^*\|_2$ stay significantly below the theoretical bound and also show an exponential decay.

**7. Conclusion.** In the present paper, we studied quasi-variational inequalities for solution-driven adaptive image denoising and non-blind image deblurring. Our general approach covers various adaptive and anisotropic types of TV regularization of first- and higher-order.

We provided theory for uniqueness and showed convergence of suitable algorithms for a broad sub-class of the considered QVIs, namely those, for which the operator corresponding to the fixed-point problem of the QVI is a contraction. Moreover, we provided convergence results, which we verified in the experimental part.

Our experiments show, that solution-driven adaptivity is able to improve the restoration results compared to its data-driven pendant.

Future work will focus on extensions to non-local regularization.

## Appendix A. Proof of Proposition 4.4.

*Claim (i):* Let $p^*_{res} \in \Pi_{\mathcal{N}^\perp(\mathrm{A})}\mathcal{D}(p^*_{res})$ be a solution to the restricted problem (4.9), i.e.,

$$\langle \nabla F(p^*_{res}), p - p^*_{res}\rangle \geq 0, \qquad \forall p \in \Pi_{\mathcal{N}^\perp(\mathrm{A})}(\mathcal{D}(p^*_{res})). \tag{A.1}$$

For any $p^* \in \mathcal{D}(p^*_{res})$ such that $p^*_{res} = \Pi_{\mathcal{N}^\perp(\mathrm{A})}p^*$, it holds that $p^* \in \mathcal{D}(p^*) = \mathcal{D}(p^*_{res})$. Note that at least one such $p^*$ exists. We show that any such $p^*$ is a solution to the unrestricted problem (3.12).

Now let $p \in \mathcal{D}(p^*_{res}) = \mathcal{D}(p^*)$ be arbitrary. We decompose $p$ into $p = p_{res} + p_{\mathcal{N}}$, where $p_{res} := \Pi_{\mathcal{N}^\perp(\mathrm{A})}(p)$, $p_{\mathcal{N}} := \Pi_{\mathcal{N}(\mathrm{A})}(p)$. Then it follows from $\mathrm{A}\,p = \mathrm{A}\,p_{res}$ and $\mathrm{A}\,p^* = \mathrm{A}\,p^*_{res}$ that

$$\begin{aligned} &\langle \nabla F(p^*), p - p^*\rangle && \text{(A.2)}\\ &= \langle \mathrm{A}^\top(\mathrm{A}\,p^* - f), p - p^*\rangle && \text{(A.3)}\\ &= \langle \mathrm{A}\,p^* - \mathrm{A}, \mathrm{A}(p - p^*)\rangle && \text{(A.4)}\\ &= \langle \mathrm{A}\,p^*_{res} - f, \mathrm{A}(p_{res} - p^*_{res})\rangle \overset{(A.1)}{\geq} 0. && \text{(A.5)} \end{aligned}$$

Thus $p^*$ is a solution of (3.12).

*Claim (ii):* Let $p^*$ be a solution to the problem (3.12). In particular, $p^* \in \mathcal{D}(p^*)$. We consider the decomposition $p^* = p^*_{res} + p^*_{\mathcal{N}}$, $p^*_{res} := \Pi_{\mathcal{N}^\perp(\mathrm{A})}(p^*)$, $p^*_{\mathcal{N}} := \Pi_{\mathcal{N}(\mathrm{A})}(p^*)$. Then, by our assumption,

$$p^*_{res} \in \Pi_{\mathcal{N}^\perp(\mathrm{A})}(\mathcal{D}(p^*)) = \Pi_{\mathcal{N}^\perp(\mathrm{A})}(\mathcal{D}(p^*_{res})).$$

It remains to show, that $p^*_{res}$ solves the restricted problem (4.9). To this end, let $p_{res} \in \Pi_{\mathcal{N}^\perp(\mathrm{A})}(\mathcal{D}(p^*_{res}))$ be arbitrary. There exists $p \in \mathcal{D}(p^*_{res})$ such that

$$p_{res} = \Pi_{\mathcal{N}^\perp(\mathrm{A})}p = p - \Pi_{\mathcal{N}(\mathrm{A})}p. \tag{A.6}$$

Then,

$$\langle \nabla F(p^*_{res}), p_{res} - p^*_{res} \rangle \tag{A.7}$$
$$= \langle \mathrm{A}\, p^*_{res} - f, \mathrm{A}(p_{res} - p^*_{res}) \rangle \tag{A.8}$$
$$= \langle \mathrm{A}\, p^* - f, \mathrm{A}(p_{res} - p^*) \rangle \tag{A.9}$$
$$\overset{(A.6)}{=} \langle \mathrm{A}\, p^* - f, \mathrm{A}(p - \Pi_{\mathcal{N}(\mathrm{A})} p - p^*) \rangle \tag{A.10}$$
$$= \langle \mathrm{A}\, p^* - f, \mathrm{A}(p - p^*) \rangle \tag{A.11}$$
$$= \langle \nabla F(p^*), p - p^* \rangle \geq 0, \tag{A.12}$$

where the last inequality holds, since $p \in \mathcal{D}(p^*)$ due to $\mathcal{D}(p^*) = \mathcal{D}(p^*_{res})$ and $p^*$ solves QVI (3.12). Thus $p^*_{res}$ is a solution to the restricted problem (4.9).

## Appendix B. Proof of Theorem 4.6.

The proof follows the proof of Thm. 6 in Nesterov's paper with $B = \mathrm{Id}$, but uses the specific form of $g(p) = \nabla F(p) = \mathrm{A}^\top(\mathrm{A}\, p - f)$. In particular we do not require $g$ to be a strongly monotone operator.

We fix two different points $v_1, v_2 \in \mathrm{im}(\mathrm{A})$. Let $\mathcal{D}_i := \tilde{\mathcal{D}}(v_i)$, $p_i \in T(v_i)$ and $g_i = \nabla F(p_i) = \mathrm{A}^\top(\mathrm{A}\, p_i - f)$.

If $\mathrm{A}(p_1 - p_2) = 0$, we immediately find

$$\| \mathrm{A}\, T(v_1) - \mathrm{A}\, T(v_2)\|_2 = \| \mathrm{A}\, p_1 - \mathrm{A}\, p_2\|_2 = 0 < \|v_1 - v_2\|_2. \tag{B.1}$$

Let us now assume $\mathrm{A}(p_1 - p_2) \neq 0$. Since $p_i$, $i = 1, 2$ solve $\arg\min_{p \in \mathcal{D}(v_i)} \frac{1}{2}\| \mathrm{A}\, p - f\|^2$, the VI

$$\langle \nabla F(p_i), q - p_i \rangle = \langle g_i, q - p_i \rangle \geq 0 \qquad \forall q \in \mathcal{D}_i \tag{B.2}$$

holds, and, for arbitrary large $\tau \geq 0$,

$$p_i = \Pi_{\mathcal{D}_i}(p_i - \tau g_i). \tag{B.3}$$

For the particular choice $q := \Pi_{\mathcal{D}_2}(p_1 - \tau g_1)$, we find from (B.3) and the variation rate condition (4.11) that

$$\|p_1 - q\|_2 = \|\Pi_{\mathcal{D}_1}(p_1 - \tau g_1) - \Pi_{\mathcal{D}_2}(p_1 - \tau g_1)\|_2 \leq \tilde{\eta}\|v_1 - v_2\|_2. \tag{B.4}$$

On the other hand, since $q$ minimizes the distance to $p_1 - \tau g_1$ within $\mathcal{D}_2$, and $p_2 \in \mathcal{D}_2$, it follows from the corresponding VI that

$$\langle q - (p_1 - \tau g_1), p_2 - q \rangle \geq 0 \tag{B.5}$$

holds. Therefore,

$$\langle q - p_1, p_2 - q\rangle \geq \tau\langle g_1, q - p_2\rangle \tag{B.6}$$
$$= \tau\langle g_1, q - p_1\rangle + \underbrace{\tau\langle g_2, p_1 - q\rangle}_{\geq 0} + \tau\langle g_2, q - p_2\rangle + \tau\langle g_1 - g_2, p_1 - p_2\rangle \tag{B.7}$$
$$\overset{(B.2), q\in\mathcal{D}_2}{\geq} \tau\langle g_1 - g_2, q - p_1\rangle + \tau\langle g_1 - g_2, p_1 - p_2\rangle \tag{B.8}$$
$$= \tau\langle \mathrm{A}^\top \mathrm{A}(p_1 - p_2), q - p_1\rangle + \tau\langle \mathrm{A}\, p_1 - \mathrm{A}\, p_2, \mathrm{A}\, p_1 - \mathrm{A}\, p_2\rangle \tag{B.9}$$
$$= \tau\langle \mathrm{A}^\top \mathrm{A}(p_1 - p_2), q - p_1\rangle + \tau\|\mathrm{A}\, p_1 - \mathrm{A}\, p_2\|_2^2. \tag{B.10}$$

Rearranging the terms, we find

$$\tau\|\mathrm{A}\, p_1 - \mathrm{A}\, p_2\|_2^2 \leq \langle q - p_1, p_2 - q\rangle + \tau\langle Op^\top \mathrm{A}(p_1 - p_2), p_1 - q\rangle \tag{B.11}$$
$$\leq \|q - p_1\|_2 \cdot \|q - p_2\|_2 + \tau\|\mathrm{A}\, p_1 - \mathrm{A}\, p_2\|_2 \cdot \|\mathrm{A}\, p_1 - \mathrm{A}\, q\|_2, \tag{B.12}$$

and, by dividing by $\tau\|\mathrm{A}\, p_1 - \mathrm{A}\, p_2\|_2 > 0$,

$$\|\mathrm{A}\, p_1 - \mathrm{A}\, p_2\|_2 \leq \frac{\|q - p_1\|_2 \cdot \|q - p_2\|_2}{\tau\|\mathrm{A}\, p_1 - \mathrm{A}\, p_2\|_2} + \|\mathrm{A}\, p_1 - \mathrm{A}\, q\|_2. \tag{B.13}$$

Since $\tau$ can be chosen arbitrarily large, we find

$$\|\mathrm{A}\, p_1 - \mathrm{A}\, p_2\|_2 \leq \|\mathrm{A}\, p_1 - \mathrm{A}\, q\|_2. \tag{B.14}$$

Moreover, by (B.4) and Assumption 4.5 (i) we have

$$\|\mathrm{A}\, p_1 - \mathrm{A}\, q\|_2 \leq \sqrt{\mu_2}\|p_1 - q\|_2 \overset{(B.4)}{\leq} \sqrt{\mu_2}\tilde{\eta}\|v_1 - v_2\|_2. \tag{B.15}$$

Combining (B.14) and (B.15), we find

$$\|\mathrm{A}\, T(v_1) - \mathrm{A}\, T(v_2)\|_2 = \|\mathrm{A}\, p_1 - \mathrm{A}\, p_2\|_2 \leq \sqrt{\mu_2}\|p_1 - q\|_2 \leq \sqrt{\mu_2}\tilde{\eta}\|v_1 - v_2\|_2. \tag{B.16}$$

## Appendix C. Lemma C.1.

The proof of Lemma 5.1 (error bounds) requires the following additional lemma.

Lemma C.1. *Let $F(p) := \frac{1}{2}\|\mathrm{A}\, p - f\|_2^2$ and $\mathcal{D}$ be a non-empty, closed and convex set. The minimizer $\overline{p}$ of the constrained problem $\min_{p\in\mathcal{D}} F(p)$ satisfies*

$$\frac{1}{2}\|\mathrm{A}\, p - \mathrm{A}\,\overline{p}\|_2^2 \leq F(p) - F(\overline{p}), \text{ for every } p \in \mathcal{D}. \tag{C.1}$$

*Proof.* The proof follows [15, Eqn. (20)-(25)]. We consider the decomposition of the primal functional $E(u) = H(u) + G(-\mathrm{A}^\top u)$ (cf. (3.6)), with

$$H(u) = \tfrac{1}{2}\|u - f\|_2^2, \tag{C.2}$$
$$G(u) = \sup_{p\in\mathcal{D}}\langle u, p\rangle. \tag{C.3}$$

For the Fenchel-dual $E^*(p)$ of $E$, we have $E^*(p) = H^*(\mathrm{A}\,p) + G^*(p)$, where

$$H^*(v) = \tfrac{1}{2}\|f - v\|_2^2 - \tfrac{1}{2}\|f\|_2^2, \tag{C.4}$$

$$G^*(p) = \begin{cases} 0 & \text{if } p \in \mathcal{D}, \\ \infty & \text{else.} \end{cases} \tag{C.5}$$

We note that for every $p \in \mathcal{D}$ we have $E^*(p) = F(p) - \frac{1}{2}\|f\|_2^2$. Therefore, it suffices to show (C.1) for $E^*$ instead of $F$. Now, let

$$I(p) := H^*(\mathrm{A}\,p) - H^*(\mathrm{A}\,\overline{p}) - \langle \mathrm{A}^\top(D_v H^*(\mathrm{A}\,\overline{p})), p - \overline{p}\rangle \tag{C.6}$$

$$J(p) := G^*(p) - G^*(\overline{p}) + \langle \mathrm{A}^\top(D_v H^*(\mathrm{A}\,\overline{p})), p - \overline{p}\rangle. \tag{C.7}$$

Then, by definition

$$I(p) + J(p) = E^*(p) - E^*(\overline{p}). \tag{C.8}$$

Since $H^*(v) = \frac{1}{2}\|v - f\|_2^2$ is strongly convex with parameter 1, i.e., $H^*(v) - H^*(v') - \langle D_v H^*(v'), v - v'\rangle \geq \frac{1}{2}\|v - v'\|_2^2$ , we find

$$H^*(\mathrm{A}\,p) - H^*(\mathrm{A}\,\overline{p}) - \langle D_v H^*(\mathrm{A}\,\overline{p}), \mathrm{A}\,p - \mathrm{A}\,\overline{p}\rangle \geq \frac{1}{2}\|\mathrm{A}\,p - \mathrm{A}\,\overline{p}\|_2^2 \tag{C.9}$$

$$\Leftrightarrow \quad H^*(\mathrm{A}\,p) - H^*(\mathrm{A}\,\overline{p}) - \langle \mathrm{A}^\top(D_v H^*(\mathrm{A}\,\overline{p})), p - \overline{p}\rangle \geq \frac{1}{2}\|\mathrm{A}\,p - \mathrm{A}\,\overline{p}\|_2^2 \tag{C.10}$$

$$\Leftrightarrow \quad I(p) \geq \|\mathrm{A}\,p - \mathrm{A}\,\overline{p}\|_2^2. \tag{C.11}$$

Now we show that $J(p) \geq 0$. Since $\overline{p}$ is the minimizer of $E^*$, we have

$$0 \in \partial(H^* \circ \mathrm{A})(\overline{p}) + \partial G^*(\overline{p}) \tag{C.12}$$

$$\Leftrightarrow \quad -\mathrm{A}^\top(D_v H^*(\mathrm{A}\,\overline{p})) \in \partial G^*(\overline{p}), \tag{C.13}$$

where $\partial H^* = \{\mathrm{A}^\top(D_v H^*(\mathrm{A}\,\overline{p}))\}$ and $\partial G^*$ are the sub-differentials (cf. [33]) of $H^*$ and $G^*$, respectively. By definition of the sub-differential,

$$-\mathrm{A}^\top(D_p H^*(\mathrm{A}\,\overline{p})) \in \partial G^*(\overline{p}) \tag{C.14}$$

$$\Leftrightarrow \quad G^*(p) \geq G^*(\overline{p}) + \langle -\mathrm{A}^\top(D_p H^*(\mathrm{A}\,\overline{p})), p - \overline{p}\rangle, \qquad \forall p, \tag{C.15}$$

$$\Leftrightarrow \quad J(p) \geq 0, \qquad \forall p. \tag{C.16}$$

Using (C.11) and (C.16) in (C.8) shows the claim. ■

## Appendix D. Proof of Prop. 5.4.

Recall that in the $k$-th outer iteration we solve the inner problem (5.1) with initial value $p^{(0)} = p^{[k]}$. The numerical solution $p^{(N)}$ of this inner problem is denoted as $p^{[k+1]}$. Thus, the required error bound (5.17) can be rewritten as

$$\|\mathrm{A}\,p^{[k+1]} - \mathrm{A}\,T(v^{[k]})\|_2 \leq \frac{\delta}{4}\|\mathrm{A}\,p^{[k]} - \mathrm{A}\,T(v^{[k]})\|_2, \tag{D.1}$$

or equivalently, with $v^{[k]} := \mathrm{A}\, p^{[k]}$,

$$\|v^{[k+1]} - \mathrm{A}\, T(v^{[k]})\|_2 \leq \frac{\delta}{4}\|v^{[k]} - \mathrm{A}\, T(v^{[k]})\|_2 = \tfrac{\delta}{4} r_k, \tag{D.2}$$

where $r_k := \|v^{[k]} - \mathrm{A}\, T(v^{[k]})\|_2 = \|\mathrm{A}\, p^{[k]} - \mathrm{A}\, T(v^{[k]})\|_2$. The proof of Prop. 5.4 follows the lines of the proof of Thm. 4.4 by Nesterov & Scrimali. Let $v^*$ be the unique fixed-point of $\mathrm{A}\, T(v)$ provided by Thm. 4.6 (2). We have

$$\begin{aligned} r_k &\geq \|v^{[k]} - v^*\|_2 - \|\mathrm{A}\, T(v^*) - \mathrm{A}\, T(v^{[k]})\|_2 \\ &\overset{\text{Thm.4.6(1)}}{\geq} \|v^{[k]} - v^*\|_2 - \tilde{\eta}\sqrt{\mu_2}\|v^{[k]} - v^*\|_2 = \delta\|v^{[k]} - v^*\|_2. \end{aligned} \tag{D.3}$$

Now we show $r_k \leq \exp(-\frac{1}{2}k) r_0$:

$$\begin{aligned} r_{k+1} &= \|v^{[k+1]} - \mathrm{A}\, T(v^{[k+1]})\|_2 && \text{(D.4)} \\ &\leq \|v^{[k+1]} - \mathrm{A}\, T(v^{[k]})\|_2 + \|\mathrm{A}\, T(v^{[k+1]} - \mathrm{A}\, T(v^{[k]})\|_2 && \text{(D.5)} \\ &\overset{\text{(D.2)}}{\leq} \tfrac{\delta}{4}\|v^{[k]} - \mathrm{A}\, T(v^{[k]})\|_2 + \tilde{\eta}\sqrt{\mu_2}\|v^{[k+1]} - v^{[k]}\|_2 && \text{(D.6)} \\ &\leq \tfrac{\delta}{4}\|v^{[k]} - \mathrm{A}\, T(v^{[k]})\|_2 + \tilde{\eta}\sqrt{\mu_2}\|v^{[k+1]} - \mathrm{A}\, T(v^{[k]})\|_2 + \tilde{\eta}\sqrt{\mu_2}\|v^{[k]} - \mathrm{A}\, T(v^{[k]})\|_2 && \text{(D.7)} \\ &\overset{\text{(D.2)}}{\leq} (\tfrac{\delta}{4} + \tfrac{\delta}{4}\tilde{\eta}\sqrt{\mu_2} + \tilde{\eta}\sqrt{\mu_2}) r_k. && \text{(D.8)} \end{aligned}$$

Inserting $\tilde{\eta}\sqrt{\mu_2} = 1 - \delta$ we derive

$$\begin{aligned} r_{k+1} &\leq (\tfrac{\delta}{4} + 1 - \delta + \tfrac{\delta}{4}(1-\delta)) r_k \\ &= (1 - \tfrac{\delta}{2} - \tfrac{\delta^2}{4}) r_k \end{aligned} \tag{D.9}$$

Applying (D.9) recursively and using $(1-s)^n \leq \exp(-sn)$, we find

$$r_k \leq \exp(-(\tfrac{\delta}{2} + \tfrac{\delta^2}{4})k) r_0. \tag{D.10}$$

Combining (D.3) and (D.10), we finally obtain

$$\delta\|v^{[k]} - v^*\|_2 \leq \exp(-(\tfrac{\delta}{2} + \tfrac{\delta^2}{4})k) r_0 \leq \exp(-\tfrac{\delta}{2}k) r_0 = \exp(-\tfrac{\delta}{2}k)\|v^{[0]} - v^*\|_2. \tag{D.11}$$

Note that by Assumption 4.5(iii), we have $\tilde{\eta}\sqrt{\mu} < 1$ and consequently $\delta = 1 - \tilde{\eta}\sqrt{\mu} > 0$. Dividing (D.11) by $\delta$ gives the claimed error estimate.

## Appendix E. Analytic Solution for Example 6.1.

Let $X = \{(a,b) \mid 0 \leq a \leq b \leq 1\}$. We calculate the unique fixed-point of $(a^0, b^0) \in X \to \arg\min_{(a,b)\in X} E(a,b)$, where

$$E(a,b) = Na^2 + \frac{N}{2}(b-1)^2 + 2\tilde{\alpha}(a^0, b^0)(b-a) \tag{E.1}$$

and

$$\tilde{\alpha}(a^0, b^0) = \max\{\alpha_0(1 - \kappa(b^0 - a^0)|), \varepsilon\}. \tag{E.2}$$

Recall that we assume $0 \leq a \leq b \leq 1$. With the assumption that $\kappa \leq 1 - \frac{\varepsilon}{\alpha_0}$, (E.2) simplifies to

$$\tilde{\alpha}(a^0, b^0) = \alpha_0(1 - \kappa(b^0 - a^0)). \tag{E.3}$$

For a fixed $\tilde{\alpha}(a^0, b^0)$, we find

$$\frac{\partial E(a,b)}{\partial a} = 2Na - 2\tilde{\alpha}(a^0, b^0), \tag{E.4}$$

$$\frac{\partial E(a,b)}{\partial b} = N(b-1) + 2\tilde{\alpha}(a^0, b^0). \tag{E.5}$$

$$\tag{E.6}$$

From $\nabla E(\bar{a}, \bar{b}) = 0$ for optimal $\bar{a}, \bar{b}$ it follows

$$\bar{a} = \frac{\tilde{\alpha}(a^0, b^0)}{N}, \qquad \bar{b} = 1 - \frac{2\tilde{\alpha}(a^0, b^0)}{N} \tag{E.7}$$

and $\tilde{\alpha}(\bar{a}, \bar{b}) = \alpha_0(1 - \kappa(1 - \frac{3\tilde{\alpha}(a^0,b^0)}{N}))$. Thus, for a fixed-point $(a^*, b^*) := (\bar{a}, \bar{b}) = (a^0, b^0)$:

$$\tilde{\alpha}(a^*, b^*) = \alpha_0 \left(1 - \kappa \left(1 - \frac{3\tilde{\alpha}(a^*, b^*)}{N}\right)\right) \tag{E.8}$$

$$\Leftrightarrow \tag{E.9}$$

$$\tilde{\alpha}(a^*, b^*) = \alpha_0 \frac{1 - \kappa}{1 - \frac{3\alpha_0\kappa}{N}}. \tag{E.10}$$

In order to guarantee $0 \leq \tilde{\alpha}(a^*, b^*) < \infty$, it suffices to show that $\kappa \leq 1$ and $\alpha_0\kappa < \frac{N}{3}$. The first condition is already covered by our assumption $\kappa \leq 1 - \frac{\varepsilon}{\alpha_0}$. We note that the second condition is weaker that the condition $\kappa\alpha_0 < \frac{1}{\mu_2} \approx \frac{1}{4}$ ($\mu_2 \approx 4$ for $d = 1$ and $n = 150$), which we found to guarantee uniqueness for the more general fixed-point problem in Example 4.8.